
\magnification=1000
\def\couleurs{oui}
\def\arXiv{oui}
  
 \def\oui{oui} 
  
\ifx\arXiv\oui
\else
 \pdfpagewidth=210truemm
 \pdfpageheight=297truemm 
\fi
  
  %
  %

  %
  %

  \catcode`@=12 

 \def\defrefnote#1{\definexref{#1}{{\the\footnotenumber}}{refnotes}}

  %
  %


\ifx\couleurs\oui
\input graphicx
 \pdfpagewidth=210truemm
 \pdfpageheight=297truemm 
 \voffset=-5mm
\fi

\input eplain.tex
\expandafter\def\expandafter\newdimen\expandafter{\newdimen}

\ifx\couleurs\oui
\beginpackages
\usepackage{color}
\endpackages 
 \pdfpagewidth=210truemm
 \pdfpageheight=297truemm 
\long\def\rge#1{{\color{red}#1}}

\definecolor{bleu-iecn}{cmyk}{.98,.13,.1,.55}

\else
\long\def\rge#1{#1}

\fi

\makeatletter
\def\numberedfootnote{%
ÊÊ\global\advance\footnotenumber by 1
ÊÊ\@eplainfootnote{{\number\footnotenumber}}%
}%
\def\makecolumns#1/#2 {\par \begingroup
ÊÊ \@columndepth = #1
ÊÊ \advance\@columndepth by -1
ÊÊ \divide \@columndepth by #2
ÊÊ \advance\@columndepth by 1
ÊÊ \@linestogoincolumn = \@columndepth
ÊÊ \@linestogo = #1
ÊÊ \currentcolumn = 1
ÊÊ \def\@endcolumnactions{%
ÊÊÊÊÊÊ\ifnum \@linestogo<2
ÊÊÊÊÊÊÊÊ \the\crtok \egroup \endgroup \par 
ÊÊÊÊÊÊ\else
ÊÊÊÊÊÊÊÊ \global\advance\@linestogo by -1
ÊÊÊÊÊÊÊÊ \ifnum\@linestogoincolumn<2
ÊÊÊÊÊÊÊÊÊÊÊÊ\global\advance\currentcolumn by 1
ÊÊÊÊÊÊÊÊÊÊÊÊ\global\@linestogoincolumn = \@columndepth
ÊÊÊÊÊÊÊÊÊÊÊÊ\the\crtok
ÊÊÊÊÊÊÊÊ \else
ÊÊÊÊÊÊÊÊÊÊÊÊ&\global\advance\@linestogoincolumn by -1
ÊÊÊÊÊÊÊÊ \fi
ÊÊÊÊÊÊ\fi
ÊÊ }%
ÊÊ \makeactive\^^M
ÊÊ \letreturn \@endcolumnactions
ÊÊ \@columnwidth = \hsize
ÊÊÊÊ \advance\@columnwidth by -\parindent
ÊÊÊÊ \divide\@columnwidth by #2
ÊÊ \penalty\abovecolumnspenalty
ÊÊ \noindent 
ÊÊ \valign\bgroup
ÊÊÊÊ &\hbox to \@columnwidth{\strut \hsize = \@columnwidth ##\hfil}\cr
}%
\makeatother

\lefteqnumbers
   \def\testd{oui}
   \def\choixlat{\ifx\numadroite\testd\righteqnumbers
            \else  \lefteqnumbers\fi}
    \choixlat

\catcode`@=\letter
\def\@eplainfootnote#1{\let\@sf\empty 
  \ifhmode\edef\@sf{\spacefactor\the\spacefactor}\/\fi
  \global\advance\hlfootlabelnumber by 1
  \hlstart@impl{foot}{\hlfootlabel}%
  \hldest@impl{footback}{\hlfootbacklabel}%
  \hbox{$^{(#1)}$}%
  \hlend@impl{foot}%
  \@sf\vfootnote{#1.}%
}%
\catcode`@=\other

  \interfootnoteskip=0pt
  \let\note=\numberedfootnote
  \everyfootnote={\eightpoint\leftskip=5truemm\rightskip5truemm}
  
  \hsize150truemm\vsize 240truemm\hoffset=5truemm

  \pretolerance=500\tolerance=1000\brokenpenalty=5000
  \parindent3mm
  
  \countdef\temps=170
  \temps=\time
  \countdef\nminutes=171{\nminutes = \time}
  \countdef\nheure=172
  \def\heure{\begingroup                     
     \temps = \time \divide\temps by 60
     \nheure = \temps                        
     \nminutes = \time
     \multiply\temps by 60
     \advance\nminutes by -\temps            
     \ifnum\nminutes<10 \toks1 = {0}%
     \else\toks1 = {}%
     \fi
     \number\nheure h\the\toks1 \number\nminutes  
  \endgroup}%

  \newcount\chstart
  \chstart=\pageno
 \headline={\ifnum\pageno=\chstart {\hfill} \else {\hss \tenrm --\ \folio\ --\hss}\fi}
  \footline={\hfill}
  \normalbaselines
  \frenchspacing
    \def\dater{\vglue-10mm\rightline{(\the\day/\the\month/\the\year)}}
  \def\dateheure{\vglue-10mm\rightline{(\the\day/\the\month/\the\year,\ \heure)}}

  \newif\ifpagetitre \pagetitretrue
\newtoks\hautpagetitre \hautpagetitre={\hfill}
\newtoks\baspagetitre \baspagetitre={\hfill}
\newtoks\auteurcourant \auteurcourant={\hfill}
\newtoks\titrecourant \titrecourant={\hfill}
\newtoks\hautpagegauche
\newtoks\hautpagedroite
\newtoks\hautpagemilieu
\hautpagemilieu={\tenrm\hfil -- \folio\ -- \hfil}
\hautpagegauche={\ifx\midfolio\oui\the\hautpagemilieu\else\tenrm\folio\hfill\the\auteurcourant\hfill\fi}
\hautpagedroite={\ifx\midfolio\oui\the\hautpagemilieu\else\hfill\the\titrecourant\hfill\tenrm\folio\fi}
\newtoks\baspagegauche \baspagegauche={\hfil}
\newtoks\baspagedroite \baspagedroite={\hfil}
\headline={\ifpagetitre\the\hautpagetitre
\else\ifodd\pageno\the\hautpagedroite\else\the\hautpagegauche\fi\fi }
\footline={\ifpagetitre\the\baspagetitre
\else\ifodd\pageno\the\baspagedroite
\else\the\baspagegauche\fi\fi \global\pagetitrefalse}

\def\pageblanche{\vfill\eject\pagetitretrue
\null\vfill\eject
\pagetitretrue
}
\def\chgtpage{\ifodd\pageno \else
\pageblanche \fi \pagetitretrue\titreun=0\footnotenumber=0}

\def\chgtpageincrtitreun{\ifodd\pageno \else
\pageblanche \fi \pagetitretrue\footnotenumber=0}

\def\majnombres{\ifodd\pageno \else
\pageblanche \fi \pagetitretrue\hautpoly\titreun=0\footnotenumber=0}

\def\hautspages#1#2{\auteurcourant={\ninepcap#1}\titrecourant={\nineit#2}}

\ifnum\chstart=\pageno \pagetitretrue\fi
  

  \def\PAR{\par}
  
  \def\leftnote#1{\vadjust{\setbox1=\vtop{\hsize 20mm\parindent=0pt\eightpoint
  \baselineskip=9pt\rightskip=4mm plus 4mm\vskip-4.7mm#1}\hbox{\kern-2cm\smash{\box1}}}}

  
  \def\raggedcenter{\leftskip=20pt plus 10em  
       \rightskip=\leftskip 
        \parfillskip=0pt 
         \spaceskip=.3333em \xspaceskip=.5em 
          \pretolerance=9999 \tolerance=9999
           \hyphenpenalty=9999 \exhyphenpenalty=9999 }
           
  \def\titrecentre#1{{\parindent0mm\raggedcenter
       \spaceskip=.6em plus .2em minus .2em\xspaceskip=.6em plus .2em minus .2em
        \tit#1\par}}
        


  \def\oui{oui}
  
\def\fontetitreun{\ifx\paradouze\oui\douzepts\gpdouze\twelvebf\textfont1=\twelveib\else
\quatorzepts\gpquatorze\fourteenbf\fi}

\def\fontetitreunl{\douzepts\textfont1=\twelveib\scriptfont1=\tenib\fourteenti}
 
 \def\fontetitredeux{\textfont1=\eleveni\ifx\paradouze\oui\onzepts\scriptfont1=\ninei\elevenit\else
                        \douzepts\twelveit\fi}
 
   \def\fontetitredeuxb{\ifx\paradouze\oui\onzepts\eleventi\gponze\textfont1=\elevenib\scriptfont1=\nineib
                         \else\douzepts\twelveti\scriptfont1=\twelveib\scriptfont1=\tenib\gpdouze\fi}
                         
\def\fontetitredeuxl{\onzepts\textfont1=\elevenbf\scriptfont1=\ninebf\twelvebf}
  
\def\fontetitretrois{\textfont0=\elevenrm\scriptfont0=\eightrm\textfont1=\eleveni
                      \scriptfont1=\eighti\scriptscriptfont1=\sixi\elevenit}
                      
\def\fontetitrequatre{\textfont0=\elevenrm\scriptfont0=\eightrm\textfont1=\eleveni
                      \scriptfont1=\eighti\scriptscriptfont1=\sixi\elevenrm}
  
  \newcount\titreun\titreun=0
  \newcount\titredeux\titredeux=0
  \newcount\titretrois\titretrois=0
  \newcount\titrequatre\titrequatre=0
  \newcount\enonce\enonce=0
  
  \def\incr#1{\global\advance#1 by 1 {\the #1}}
  \def\avance#1{\global\advance#1 by 1}
  \def\init#1{\global#1=0}
  
  \long\def\Indentation#1#2{\setbox10=\hbox{\fontetitreun#1}
  	                    \ifdim\wd10 < 4mm
                         \setbox10=\hbox to 4mm{\box10\hfill}
                       \else\ifdim\wd10 < 6mm
                         \setbox10=\hbox to 6mm{\box10\hfill}
  	                    \else\ifdim\wd10 < 8mm
                         \setbox10=\hbox to 8mm{\box10\hfill}
                       \else\ifdim\wd10 < 12mm
                         \setbox10=\hbox to 12mm{\box10\hfill}
                       \fi\fi\fi\fi
                       \dimen10=\hsize
                       \advance \dimen10 by -\wd10
                       \noindent \box10 %
                       \ignorespaces
                       \hbox{\vtop{\hsize=\dimen10\raggedright\noindent\fontetitreun#2}}}

  \long\def\paraun#1{\removelastskip\par\medskip\goodbreak\vskip0pt plus.01\vsize\penalty-100
                \vskip0pt plus-.01\vsize
  	              \init{\titredeux}\ifnum\optionparag=1{\init\eqnumber\init\enonce}\else{}\fi
                  \goodbreak{\fontetitreun
  	                \Indentation{\incr{\titreun}.\ }{\fontetitreun #1\par}}\nobreak\medskip}

 %
 %
 \long\def\paraunc#1{\removelastskip\par\bigskip\goodbreak\vskip0pt plus.01\vsize\penalty-100
                \vskip0pt plus-.01\vsize
  	              \init{\titredeux}
                 \ifnum\optionparag=1{\init{\eqnumber}\init\enonce}\else{}\fi
                  \goodbreak
  	                {\parindent0mm\raggedcenter\fontetitreun\incr{\titreun}.\ 
                     \fontetitreun #1\par}\nobreak\medskip}
                     
\newtoks\titreunl
\titreunl={\ifnum\titreun=1{I}\fi%
\ifnum\titreun=2{II}\fi%
\ifnum\titreun=3{III}\fi%
\ifnum\titreun=4{IV}\fi%
\ifnum\titreun=5{V}\fi%
\ifnum\titreun=6{VI}\fi%
\ifnum\titreun=7{VII}\fi%
\ifnum\titreun=8{VIII}\fi%
\ifnum\titreun=9{IX}\fi%
\ifnum\titreun=10{X}\fi%
\ifnum\titreun=11{XI}\fi%
\ifnum\titreun=12{XII}\fi%
\ifnum\titreun=13{XIII}\fi%
}
\long\def\paraunl#1{\removelastskip\par\bigskip\bigskip\goodbreak\vskip0pt plus.01\vsize\penalty-100
                \vskip0pt plus-.01\vsize
  	              \init{\titredeux}\ifnum\optionparag=1{\init\eqnumber\init\enonce}\else{}\fi
                  \goodbreak{\fontetitreunl
  	                \Indentation{\global\advance\titreun by 1{\the\titreunl}.\ }{\fontetitreunl #1\par}}\nobreak\smallskip}

  
  \long\def\paradeux#1{\init{\titretrois}\vskip0pt plus.01\vsize\penalty-10
                \vskip0pt plus-.01\vsize\ifx \elie\oui\medskip\ifnum\titredeux>0\medskip\fi\fi
                 \Indentation{\fontetitredeux\the\titreun${\cdot}$\incr{\titredeux}.}
                              {\fontetitredeux\textfont1=\eleveni#1}\nobreak\par }
  
  \long\def\paradeuxb#1{\init{\titretrois}\vskip0pt plus.001\vsize\penalty-10
                \vskip0pt plus-.01\vsize{\ifx \elie\oui\medskip\ifnum\titredeux>0\medskip\fi\fi
                  \Indentation
  {\fontetitredeuxb\the\titreun${\cdot}$\incr{\titredeux}.}{ \fontetitredeuxb#1}}\nobreak
\smallskip}

\newtoks\titredeuxl
\titredeuxl={\ifnum\titredeux=1{A}\fi%
\ifnum\titredeux=2{B}\fi%
\ifnum\titredeux=3{C}\fi%
\ifnum\titredeux=4{D}\fi%
\ifnum\titredeux=5{E}\fi%
\ifnum\titredeux=6{F}\fi%
\ifnum\titredeux=7{G}\fi%
\ifnum\titredeux=8{H}\fi%
\ifnum\titredeux=9{I}\fi%
\ifnum\titredeux=10{J}\fi%
\ifnum\titredeux=11{K}\fi%
\ifnum\titredeux=12{L}\fi%
\ifnum\titredeux=13{M}\fi%
}
 \long\def\paradeuxl#1{\init{\titretrois}\vskip0pt plus.001\vsize\penalty-10
                \vskip0pt plus-.01
                \vsize \bigskip%
                  \Indentation
     {\fontetitredeuxl\global\advance\titredeux by 1
  \quad \the\titreunl${\cdot}$\the\titredeuxl.}{ \fontetitredeuxl#1}
  \removelastskip\nobreak\smallskip}
  

  \long\def\paratrois#1{\init{\titrequatre}\ifdim\lastskip<\smallskipamount
                \removelastskip\smallskip\fi
                 \vskip0pt plus.01\vsize\penalty-10
                  \vskip0pt
plus-.01\vsize{\ifx \elie\oui\ifnum\titretrois>0\medskip\fi\fi
\Indentation{\fontetitretrois\the\titreun${\cdot}$\the\titredeux${\cdot}$\incr{\titretrois}.\ }
  {\hskip0mm\baselineskip=14pt\fontetitretrois#1}\nobreak\smallskip}}
  
  
  \long\def\paratroisl#1{\init{\titrequatre}\ifdim\lastskip<\smallskipamount
                \removelastskip\fi
                 \vskip0pt plus.01\vsize\penalty-10
                  \vskip0pt
plus-.01\vsize\ifx \elie\oui\bigskip
\fi
\Indentation{\fontetitretrois\quad \quad \the\titreunl{${\cdot}$}\the\titredeuxl${\cdot}$\incr{\titretrois}.\ }
  {\hskip0mm\fontetitretrois#1}\nobreak\smallskip}


  \long\def\paraquatre#1{\ifdim\lastskip<\smallskipamount
                \removelastskip\smallskip\fi
                 \vskip0pt plus.01\vsize\penalty-10
                  \vskip0pt
                  plus-.01\vsize\par
 
\Indentation{\fontetitrequatre \the\titreun{${\cdot}$}\the\titredeux${\cdot}$\the\titretrois${\cdot}$\incr{\titrequatre}.\ }
{\hskip0mm\fontetitrequatre#1}\nobreak\smallskip}


\newtoks\titrequatrel
\titrequatrel={\ifnum\titrequatre=1{a}\fi%
\ifnum\titrequatre=2{b}\fi%
\ifnum\titrequatre=3{c}\fi%
\ifnum\titrequatre=4{d}\fi%
\ifnum\titrequatre=5{e}\fi%
\ifnum\titrequatre=6{f}\fi%
\ifnum\titrequatre=7{g}\fi%
\ifnum\titrequatre=8{h}\fi%
\ifnum\titrequatre=9{i}\fi%
\ifnum\titrequatre=10{j}\fi%
\ifnum\titrequatre=11{k}\fi%
\ifnum\titrequatre=12{l}\fi%
\ifnum\titrequatre=13{m}\fi%
}
\long\def\paraquatrel#1{\ifdim\lastskip<\smallskipamount
                \removelastskip\smallskip\fi
                 \vskip0pt plus.01\vsize\penalty-10
                  \vskip0pt
                  plus-.01\vsize{\bigskip
\Indentation{\global\advance\titrequatre by 1
\fontetitrequatre\quad \quad \quad \the\titreunl${\cdot}$\the\titredeuxl${\cdot}$\the\titretrois${\cdot}$\the\titrequatrel.\ }
{\hskip0mm\fontetitrequatre#1}\nobreak\smallskip}}

\ifx\optionkeys\oui
\def\drefun#1{\definexref{¤#1}{{\the\titreun}}{}} 
\def\drefdeux#1{\definexref{¤#1}{{\the\titreun}.{\the\titredeux}}{}}
\def\dreftrois#1{\definexref{¤#1}{{\the\titreun}.{\the\titredeux}.{\the\titretrois}}{}}
\else
\def\drefun#1{\definexref{prg#1}{{\the\titreun}}{}} 
\def\drefdeux#1{\definexref{prg#1}{{\the\titreun}.{\the\titredeux}}{}}
\def\dreftrois#1{\definexref{prg#1}{{\the\titreun}.{\the\titredeux}.{\the\titretrois}}{}}
\fi

%


  \long\def\propdeux#1#2#3#4{%
       \avance{\enonce}
       \leavevmode\edef\temp{#2}%
         \ifx\temp\empty 
          \else
           \definexref{#2}{#1~{\the\titreun.\the\enonce}}{enonces}
            \definexref{s#2}{{\the\titreun.\the\enonce}}{enonces}
             \fi
\smallskip
      \noindent{\bf#1\ {\bf\the\titreun.\the\enonce{#3}.}\enspace}{\sl#4\par}%
      \ifdim\lastskip<\medskipamount \removelastskip\penalty55\par \fi
   }

  \long\def\propun#1#2#3#4{%
      \avance{\enonce}
       \leavevmode\edef\temp{#2}%
        \ifx\temp\empty 
          \else
           \definexref{#2}{#1~{\the\enonce}}{enonces}
            \definexref{{s#2}}{{\the\enonce}}{enonces}
             \fi
   \par 
     \noindent{\bf#1\ {\bf\the\enonce{#3}.}\enspace}{\sl#4\par}%
     \ifdim\lastskip<\medskipamount \removelastskip\penalty55\medskip\fi
  }
  
  \long\def\prop#1#2#3#4{\ifnum\optionparag=1
                          \propdeux{#1}{#2}{\textfont1=\elevenib#3}{#4} \else\propun{#1}{#2}{\textfont1=\elevenib#3}{#4}\fi}

  \long\def\propt#1#2#3{\ifx\tpf\oui \prop{Th\'eo\-r\`eme}{#1}{#2}{#3}\par
                       \else\prop{Theorem}{#1}{#2}{#3}\par\fi}
  \long\def\Propt#1#2{\propt{#1}{}{#2}}
  \long\def\propl#1#2#3{\ifx\tpf\oui\prop{Lem\-me}{#1}{#2}{#3}\par
                         \else\prop{Lemma}{#1}{#2}{#3}\par\fi}
  \long\def\Propl#1#2{\propl{#1}{}{#2}}
  \long\def\propc#1#2#3{\ifx\tpf\oui\prop{Corol\-laire}{#1}{#2}{#3}\par
                         \else\prop{Corollary}{#1}{#2}{#3}\par\fi}
  \long\def\Propc#1#2{\propc{#1}{}{#2}}
  \long\def\propp#1#2#3{\prop{Pro\-po\-si\-tion}{#1}{#2}{#3}\par}
  \long\def\Propp#1#2{\propp{#1}{}{#2}} 
  \long\def\propd#1#2#3{\ifx\tpf\oui\prop{D\'efi\-nition}{#1}{#2}{#3}\par
                       \else\prop{Definition}{#1}{#2}{#3}\par\fi} 
  
  \long\def\proptd#1#2#3{\ifx\tpf\oui\prop{Th\'eor\`eme et d\'efi\-nition}{#1}{#2}{#3}\par
                       \else\prop{Theorem and definition}{#1}{#2}{#3}\par\fi}


  
  \newcount\optionparag\optionparag=1
  
  \long\def\section#1#2{\ifnum\optionparag=1 \paraun{#2} 
                        \else\goodbreak{\fontetitreun
  	                \Indentation{#1.\ }{#2}}\nobreak\smallskip\fi}

  \def\eqconstruct#1{\ifnum\optionparag=1{\the\titreun\hbox{$\cdot$}#1}\else{#1}\fi}

  
  
  \def\numref{oui}  
  
  \newcount\mesref\mesref=0 
  \def\defbib#1{\ifx\numref\oui\global\advance\mesref by 1 \definexref{#1}{{\the
                 \mesref}}{}\else\definexref{#1}{#1}{}\fi}
  \def\bibtem#1{\defbib{#1}\item{\citer{#1}}}
  \def\citer#1{[\ref{#1}]}
  \def\citeplus#1#2{[\ref{#1}; #2]}

  
  \font\seventeenmsa=msam10 at 17pt    
  \font\fourteenmsa=msam10 at 14pt
  \font\twelvemsa=msam10 at 12pt
  \font\tenmsa=msam10                 
  \font\ninemsa=msam10 at 9pt 
  \font\eightmsa=msam10 at 8pt 
  \font\sevenmsa=msam7 
  \font\sixmsa=msam10 at 6pt
  \font\fivemsa=msam5
  \newfam\msafam\textfont\msafam=\tenmsa\scriptfont\msafam=\sevenmsa\scriptscriptfont\msafam=\fivemsa
  
  \font\seventeenbb=msbm10 at 17pt     
  \font\fourteenbb=msbm10 at 14pt
  \font\twelvebb=msbm10 at 12pt
  \font\tenbb=msbm10                   
  \font\ninebb=msbm10 at 9pt 
  \font\eightbb=msbm10 at 8pt 
  \font\sevenbb=msbm7 
  \font\sixbb=msbm10 at 6pt
  \font\fivebb=msbm5 
  \newfam\bbfam\textfont\bbfam=\tenbb\scriptfont\bbfam=\sevenbb\scriptscriptfont\bbfam=\fivebb
  \def\bb{\fam\bbfam\tenbb}%

  \font\seventeenscaln=eusm10 at 17pt   
  \font\twelvescaln=eusm10 at 12pt
  \font\tenscaln=eusm10                
  \font\ninescaln=eusm10 scaled 900
  \font\eightscaln=eusm10 scaled 800
  \font\sevenscaln=eusm10 scaled 700
  \font\sixscaln=eusm10 scaled 600
   
  \newfam\scalnfam\textfont\scalnfam=\tenscaln\scriptfont\scalnfam=\sevenscaln\scriptscriptfont\scalnfam=\sixscaln
  \def\scaln{\fam\scalnfam\tenscaln}%
  \def\scal{\scaln}
  
  \font\tenscalb=eusb10                

  \font\sevenscalb=eusb10 scaled 700

  \newfam\scalbfam\textfont\scalbfam=\tenscalb\scriptfont\scalbfam=\sevenscalb
  %
  
  %
  %
  \font\fourteenrm=cmr12 scaled 1200
  \font\elevenrm=cmr10 at 11pt
  \font\twelverm=cmr12
  \font\ninerm=cmr9
  \font\eightrm=cmr8      
  \font\sevenrm=cmr7
  \font\sixrm=cmr6

  \font\seventeenpcap=cmcsc10 at 17pt
  \font\tenpcap=cmcsc10                        
  \font\ninepcap=cmcsc9
  \font\eightpcap=cmcsc8
  \font\sevenpcap=cmcsc10 scaled 700
  
  \newfam\pcapfam\textfont\pcapfam=\tenpcap\scriptfont\pcapfam=\sevenpcap
  \def\pcap{\fam\pcapfam\tenpcap}
  
  \font\seventeenrm=cmbx12 scaled 1400

  \font\fourteenbf=cmbx10 scaled 1400
  
  \font\twelvebf=cmbx12
  \font\elevenbf=cmbx10 at 11pt
  \font\ninebf=cmbx9  
  \font\eightbf=cmbx8
  \font\sixbf=cmbx6
  
  \font\tengot=eufm10                           
   
  \font\eightgot=eufm10 at 8truept 
  \font\sevengot=eufm7 
  \font\sixgot=eufm10 at 6 truept 
   
  \newfam\gotfam
  \textfont\gotfam=\tengot\scriptfont\gotfam=\sevengot\scriptscriptfont\gotfam=\sixgot
  \def\got{\fam\gotfam\tengot}%

  
  \def\tit{%
  \textfont0=\seventeenrm\scriptfont0=\tenrm\def\rm{\fam0\seventeenrm}%
  \textfont1=\seventeenib\scriptfont1=\twelveib%
  \textfont2=\seventeensy\scriptfont2=\twelvesy\scriptscriptfont2=\ninesy
  \textfont3=\seventeenex
  \textfont\itfam=\seventeenti
  \def\it{\fam\itfam\seventeenti}%
  \textfont\bbfam=\seventeenbb \scriptfont\bbfam=\twelvebb
  \def\bb{\fam\bbfam\seventeenbb}%
  \textfont\msafam=\seventeenmsa\scriptfont\msafam=\twelvemsa
  \textfont\scalnfam=\seventeenscaln
  \def\pcap{\fam\pcapfam\seventeenpcap}
  \normalbaselineskip=25pt\normalbaselines\rm}

  \font\seventeenti=cmbxti10 scaled 1680
  
  \font\fourteenti=cmbxti10 at 14pt
  
  \font\twelveti=cmbxti10 scaled 1200
  \font\eleventi=cmbxti10 at 11pt

  %
  %
  \font\twelveit=cmti12	
  \font\elevenit=cmti10 scaled 1100
  \font\nineit=cmti9
  \font\eightit=cmti8
  \font\sevenit=cmti7

  %
  %
 
 \font\seventeenib=cmmib10 scaled 1680
  \font\fourteenib=cmmib10 scaled 1400
  \font\twelveib=cmmib10 scaled 1200
  \font\elevenib=cmmib10 scaled 1100
  \font\tenib=cmmib10
\font\eightib=cmmib10 scaled 800
  \font\nineib=cmmib10 scaled 900
\font\sevenib=cmmib10 scaled 700
\font\sixib=cmmib10 scaled 600
\font\fiveib=cmmib10 scaled 500

\ifx\ITAN\oui
\else
\innernewfam\cmmibfam
\textfont\cmmibfam=\tenib
\scriptfont\cmmibfam=\sevenib
\scriptscriptfont\cmmibfam=\fiveib
\def\ib{\fam\cmmibfam\tenib}
\fi

  %
  %
  \font\twelvei=cmmi10 scaled 1200
  \font\eleveni=cmmi10 scaled 1100
  \font\ninei=cmmi9
  \font\eighti=cmmi8 
  \font\seveni=cmmi7 			                
  \font\sixi=cmmi6
  
  \font\ninesl=cmsl9                    
  \font\eightsl=cmsl8 
  \font\sevensl=cmsl10 at 7pt

  \font\ninett=cmtt9                    
  \font\eighttt=cmtt8
  \font\seventt=cmtt10 scaled 700

  \font\seventeensy=cmsy10 scaled 1680    
  \font\fourteensy=cmsy10 scaled 1400
  \font\twelvesy=cmsy10 scaled 1176
  
  \font\ninesy=cmsy9                      
  \font\eightsy=cmsy8
  \font\sixsy=cmsy6
  \font\seventeenex=cmex10 at 17pt
  \font\fourteenex=cmex10 at 14pt
  \font\twelveex=cmex10 at 12pt
  \font\nineex=cmex10 at 9pt
  \font\eightex=cmex10 at 8pt
  \font\sevenex=cmex10 at 7pt
  \font\sixex=cmex10 at 6pt
  \font\fiveex=cmex10 at 5pt
  
   
  \font\fourteengp=cmmi10 at 14pt
  
  \font\twelvegp=cmmib10 at 12pt
  \font\elevengp=cmmib10 at 11pt
  \font\tengp=cmmib10                          
  \font\ninegp=cmmib10 at 9pt 
  \font\eightgp=cmmib8 
  \font\sevengp=cmmib7 
  \font\sixgp=cmmib6


\def\gp{\textfont0=\tenbf\scriptfont0=\sevenbf\scriptscriptfont0=\fivebf
\textfont1=\tengp\scriptfont1=\sevengp\scriptscriptfont1=\fivegp}
  \def\gponze{\textfont0=\elevenbf\scriptfont0=\eightbf\scriptscriptfont0=\sixbf
           \textfont1=\elevengp\scriptfont1=\eightgp\scriptscriptfont1=\sixgp}
  \def\gpdouze{\textfont0=\twelvebf\scriptfont0=\tenbf\scriptscriptfont0=\ninebf
           \textfont1=\twelvegp\scriptfont1=\tengp\scriptscriptfont1=\ninegp}        
  
 \def\gpquatorze{\textfont0=\fourteenbf\scriptfont0=\twelvebf\scriptscriptfont0=\elevenbf
           \textfont1=\fourteengp\scriptfont1=\twelvegp\scriptscriptfont1=\elevengp}

  
  \expandafter\chardef\csname pre amssym.def at\endcsname=\the\catcode`\@
  \catcode`\@=11
  \def\undefine#1{\let#1\undefined}
  \def\newsymbol#1#2#3#4#5{\let\next@\relax
   \ifnum#2=\@ne\let\next@\msafam@\else
   \ifnum#2=\tw@\let\next@\bbfam@\fi\fi
   \mathchardef#1="#3\next@#4#5}
  \def\mathhexbox@#1#2#3{\relax
   \ifmmode\mathpalette{}{\m@th\mathchar"#1#2#3}%
   \else\leavevmode\hbox{$\m@th\mathchar"#1#2#3$}\fi}
  \def\hexnumber@#1{\ifcase#1 0\or 1\or 2\or 3\or 4\or 5\or 6\or 7\or 8\or
   9\or A\or B\or C\or D\or E\or F\fi}
  
  \def\setboxz@h{\setbox\z@\hbox}
  \def\wdz@{\wd\z@}
  \def\boxz@{\box\z@}
  
  \edef\msafam@{\hexnumber@\msafam}
  \mathchardef\dabar@"0\msafam@39
  
  \edef\bbfam@{\hexnumber@\bbfam}
  \def\widehat#1{\setboxz@h{$\m@th#1$}%
   \ifdim\wdz@>\tw@ em\mathaccent"0\bbfam@5B{#1}%
   \else\mathaccent"0362{#1}\fi}
  \def\widetilde#1{\setboxz@h{$\m@th#1$}%
   \ifdim\wdz@>\tw@ em\mathaccent"0\bbfam@5D{#1}%
   \else\mathaccent"0365{#1}\fi}
  \newsymbol\leqq 1335          
  \newsymbol\leqslant 1336
  \newsymbol\lessgtr 1337       
  \newsymbol\backprime 1038     
  \newsymbol\risingdotseq 133A  
  \newsymbol\fallingdotseq 133B 
  \newsymbol\succcurlyeq 133C   
  \newsymbol\geqq 133D          
  \newsymbol\geqslant 133E
  \newsymbol\nmid 232D
  \newsymbol\nexists 2040
  \newsymbol\smallsetminus 2272
  \newsymbol\varnothing 203F 
  \catcode`\@=\active

  \catcode`\@=11
  \newcount\typofr\typofr=1
  
  \catcode`\;=\active
  \def;{\ifnum\typofr=1\relax\ifhmode\ifdim\lastskip>\z@\unskip\fi
     \kern.2em\fi\string;\else\string;\fi}
  
  \catcode`\:=\active
  \def:{\ifnum\typofr=1\relax\ifhmode\ifdim\lastskip>\z@\unskip\fi
  \penalty\@M\ \fi\string:\else\string:\fi}
  
  \catcode`\!=\active
  \def!{\ifnum\typofr=1\relax\ifhmode\ifdim\lastskip>\z@\unskip\fi
     \kern.2em\fi\string!\else\string!\fi}
  
  \catcode`\?=\active
  \def?{\ifnum\typofr=1\relax\ifhmode\ifdim\lastskip>\z@\unskip\fi
     \kern.2em\fi\string?\else\string?\fi}

  \def\francais{\typofr=1\def\tpf{oui}}
  \def\anglais{\typofr=2\def\tpf{non}\def\english{oui}}
  \def\oui{oui}
  \francais
  
  \catcode`\@=12
  

%
\def\raye #1|{\leavevmode\setbox1=\hbox{#1}%
\raise .5pt\hbox to \wd1{\xleaders\hbox{\rge{$ \sct / $}%
\kern 1pt}\hfill\kern -1pt }\kern -\wd1 \unhbox1\relax }

\def\barre #1|{\leavevmode\setbox1=\hbox{#1}%
\rlap{\color{red}\vrule height 2.4pt depth -1.2pt width \wd1}\color{black} \unhbox1\relax}

  

  
  \def\og{\leavevmode\raise.24ex\hbox{$\scriptscriptstyle\langle\!\langle\>$}}    
  \def\fg{\leavevmode\raise.24ex\hbox{$\scriptscriptstyle\>\rangle\!\rangle$}}    

  \def\d{\,{\rm d}}
  \def\dd{{\rm d}}

  \def\z{{\bb Z}}
  \def\r{{\bb R}}
  
  \def\N{{\bb N}}

  \def\PP{{\bb P}}

  \def\B{{\scal B}}
  
  \def\D{{\scal D}}

  \def\L{{\scal L}}
  \def\M{{\scal M}}
  
  \def\O{{\scal O}}
  \def\P{{\scaln P}}

  \def\s{{\cal S}}
  \def\T{{\scal T}}
  
  \def\W{{\scal W}}

  \def\frac#1#2{{#1\over #2}}
  \def\di#1#2{\sct#1\atop{\sct#2}}

  \def\qedbox{$\rlap{$\sqcap$}\sqcup$}           
  \def\qed{\nobreak\hfill\penalty250 \hbox{}\nobreak\hfill\qedbox\par }

  \def\plaf#1{\left\lceil#1\right\rceil}
  
  \def\¤{\S\thinspace}

  \def\¥{$\bullet$ }
  
  
  \def\e{{\rm e}}

  \def\epsilon{\varepsilon}

  \def\phi{\varphi}
  \def\theta{\vartheta}
  \def\rho{\varrho}
  \def\dm{{\textstyle{1\over 2}}}

  \def\sct{\scriptstyle}
  \def\pf{\noi{\it Proof. }}
  \def\nid{\ifnum\typofr=1\par\noindent{\it D\'emonstration. }\else\pf\fi}
  \def\noi{\noindent}
  \def\rem{\ifnum\typofr=1\noi{\it Remarque.}\ \else\noi{\it Remark.}\ \fi}
  \def\rems{\ifnum\typofr=1\noi{\it Remarques.}\ \else\noi{\it Remarks.}\ \fi}

  \def\sset{\smallsetminus}

  \def\1{{\bf 1}}
  \def\|{\Vert}

  \def\leq{\leqslant}
  \def\geq{\geqslant}

  \def\ie{{i.e.\ }}
  \def\eg{{e.g.}}
  \newsymbol\subsetneqq 2324
  \newsymbol\subsetneq 2328

  \def\fl#1{\left\lfloor #1 \right\rfloor}
  \def\plaf#1{\left\lceil #1 \right\rceil}

  \def\log{\mathop{\rm log}\nolimits}




  \def\pmb#1{\setbox0=\hbox{#1}%
  \kern-.025em\copy0\kern-\wd0\kern.05em\copy0\kern-\wd0\kern-.025em\raise .0433em\box0 }

  
  \skewchar\eighti='177 \skewchar\sixi='177
  \skewchar\eightsy='60 \skewchar\sixsy='60
  
  \def\eightpoint{%
  \textfont0=\eightrm\scriptfont0=\sixrm\scriptscriptfont0=\fiverm
  \def\rm{\fam0\eightrm}%
  \textfont1=\eighti\scriptfont1=\sixi
  \scriptscriptfont1=\fivei\def\oldstyle{\fam1\seveni}%
  \textfont2=\eightsy\scriptfont2=\sixsy\scriptscriptfont2=\fivesy
  \textfont3=\eightex\scriptfont3=\sixex
  \textfont\itfam=\eightit
  \def\it{\fam\itfam\eightit}%
  \textfont\slfam=\eightsl
  \def\sl{\fam\slfam\eightsl}%
  \textfont\bbfam=\eightbb \scriptfont\bbfam=\sixbb\scriptscriptfont\bbfam=\fivebb
  \def\bb{\fam\bbfam\eightbb}%
  \textfont\msafam=\eightmsa\scriptfont\msafam=\sixmsa
  \textfont\scalnfam=\eightscaln
  \def\scaln{\fam\scalnfam\eightscaln}
  \textfont\ttfam=\eighttt
  \def\tt{\fam\ttfam\eighttt}%
\textfont\gotfam=\eightgot
  \textfont\bffam=\eightbf\scriptfont\bffam=\sixbf\scriptscriptfont\bffam=\fivebf
  \def\bf{\fam\bffam\eightbf}%
  \ifx\ITAN\oui\else\textfont\cmmibfam=\eightib
       \scriptfont\cmmibfam=\sixib
        \scriptscriptfont\cmmibfam=\fiveib
         \def\ib{\fam\cmmibfam\eightib}
   \fi
  \textfont\pcapfam=\eightpcap
  \def\pcap{\fam\pcapfam\eightpcap}
  \abovedisplayskip=2pt plus2pt minus 2pt
  \belowdisplayskip=2pt plus1pt minus 2pt
  \abovedisplayshortskip= 1pt plus 2pt minus 1pt
  \belowdisplayshortskip= 1pt plus 2pt minus 1pt
  \smallskipamount=2pt plus 1pt minus 2pt
  \medskipamount=3pt plus 2pt minus 2pt
  \bigskipamount=7pt plus 3pt minus 3pt
  \setbox\strutbox=\hbox{\vrule height 5pt depth 2pt width 0pt}%
  \normalbaselineskip=9pt\normalbaselines\rm}

  \def\({\left(}
  \def\){\right)}
  
  \def\footnoterule{\kern -2pt\hrule width 7truecm\kern 2.4pt}
  
  \def\xnotedef#1{\definexref{#1}{\noexpand\number\footnotenumber}{Note}}%

  
  
  \def\ninepoint{%
  \textfont0=\ninerm\scriptfont0=\sixrm\scriptscriptfont0=\fiverm
  \def\rm{\fam0\ninerm}%
  \textfont1=\ninei\scriptfont1=\sixi
  \scriptscriptfont1=\fivei\def\oldstyle{\fam1\ninei}%
  \textfont2=\ninesy\scriptfont2=\sixsy\scriptscriptfont2=\fivesy
  \textfont3=\nineex\scriptfont3=\sixex
  \textfont\itfam=\nineit
  \def\it{\fam\itfam\nineit}%
  \textfont\slfam=\ninesl
  \def\sl{\fam\slfam\ninesl}%
  \textfont\bbfam=\ninebb\scriptfont\bbfam=\sixbb\scriptscriptfont\bbfam=\fivebb
  \def\bb{\fam\bbfam\ninebb}%
  \textfont\msafam=\ninemsa\scriptfont\msafam=\sixmsa\scriptscriptfont\msafam=\fivemsa
  \textfont\scalnfam=\ninescaln
  \def\scaln{\fam\scalnfam\ninescaln}
  \textfont\ttfam=\ninett
  \def\tt{\fam\ttfam\ninett}%
  \textfont\bffam=\ninebf\scriptfont\bffam=\sixbf\scriptscriptfont\bffam=\fivebf
  \def\bf{\fam\bffam\ninebf}%
  \abovedisplayskip=3pt plus2pt minus 2pt
  \belowdisplayskip=3pt plus1pt minus 2pt
  \abovedisplayshortskip= 2pt plus 2pt minus 1pt
  \belowdisplayshortskip= 2pt plus 2pt minus 1pt
  \smallskipamount=2pt plus 1pt minus 2pt
  \medskipamount=3pt plus 2pt minus 2pt
  \bigskipamount=7pt plus 3pt minus 3pt
  \setbox\strutbox=\hbox{\vrule height 5pt depth 2pt width 0pt}%
  \normalbaselineskip=11pt plus.3pt minus.2pt\normalbaselines\rm}

  \def\sevenpoint{%
  \textfont0=\sevenrm\scriptfont0=\sixrm\scriptscriptfont0=\fiverm
  \def\rm{\fam0\sevenrm}%
  \textfont1=\seveni\scriptfont1=\sixi
  \scriptscriptfont1=\fivei\def\oldstyle{\fam1\seveni}%
  \textfont2=\sevensy\scriptfont2=\sixsy\scriptscriptfont2=\fivesy
  \textfont3=\sevenex\scriptfont3=\fiveex
  \textfont\itfam=\sevenit
  \def\it{\fam\itfam\sevenit}%
  \textfont\slfam=\sevensl
  \def\sl{\fam\slfam\sevensl}%
  \textfont\bbfam=\sevenbb \scriptfont\bbfam=\sixbb\scriptscriptfont\bbfam=\fivebb
  \def\bb{\fam\bbfam\sevenbb}%
  \textfont\msafam=\sevenmsa\scriptfont\msafam=\sixmsa
  \textfont\scalnfam=\sevenscaln
  \def\scaln{\fam\scalnfam\sevenscaln}
  \textfont\bffam=\sevenbf\scriptfont\bffam=\sixbf\scriptscriptfont\bffam=\fivebf
  \def\bf{\fam\bffam\sevenbf}%
  \textfont\ttfam=\seventt
  \abovedisplayskip=2pt plus2pt minus 2pt
  \belowdisplayskip=2pt plus1pt minus 2pt
  \abovedisplayshortskip= 1pt plus 2pt minus 1pt
  \belowdisplayshortskip= 1pt plus 2pt minus 1pt
  \smallskipamount=2pt plus 1pt minus 2pt
  \medskipamount=3pt plus 2pt minus 2pt
  \bigskipamount=7pt plus 3pt minus 3pt
  \setbox\strutbox=\hbox{\vrule height 5pt depth 2pt width 0pt}%
  \normalbaselineskip=9pt\normalbaselines\rm}

 \def\onzepts{%
 \textfont0=\elevenrm\scriptfont0=\ninerm
 \textfont1=\eleveni\scriptfont1=\ninei
}

\def\douzepts{%
  \textfont0=\twelverm\scriptfont0=\tenrm\def\rm{\fam0\twelverm}%
  \textfont1=\twelvei\scriptfont1=\teni%
  \textfont2=\twelvesy\scriptfont2=\tensy\scriptscriptfont2=\eightsy
  \textfont3=\twelveex
  \textfont\itfam=\twelveti
  \def\it{\fam\itfam\twelveti}%
  \textfont\bffam=\twelvebf\scriptfont\bffam=\tenbf\scriptscriptfont\bffam=\eightbf
  \def\bf{\fam\bffam\twelvebf}%
  \textfont\bbfam=\twelvebb \scriptfont\bbfam=\tenbb
  \def\bb{\fam\bbfam\twelvebb}%
  \textfont\msafam=\twelvemsa\scriptfont\msafam=\tenmsa
  \textfont\scalnfam=\twelvescaln
  \normalbaselineskip=15pt\normalbaselines\rm}

\def\quatorzepts{%
  \textfont0=\fourteenrm\scriptfont0=\twelverm\def\rm{\fam0\fourteenrm}%
  \textfont1=\fourteenib\scriptfont1=\twelveib%
  \textfont2=\fourteensy\scriptfont2=\twelvesy\scriptscriptfont2=\tensy
  \textfont3=\fourteenex
  \textfont\itfam=\fourteenti
  \def\it{\fam\itfam\fourteenti}%
  \textfont\bffam=\fourteenbf\scriptfont\bffam=\twelvebf\scriptscriptfont\bffam=\tenbf
  \def\bf{\fam\bffam\fourteenbf}%
  \textfont\bbfam=\fourteenbb \scriptfont\bbfam=\twelvebb
  \def\bb{\fam\bbfam\fourteenbb}%
  \textfont\msafam=\fourteenmsa\scriptfont\msafam=\twelvemsa
  \textfont\scalnfam=\twelvescaln
  \normalbaselineskip=18pt\normalbaselines\rm}


\def\AA{{\it Acta Arith.}}

\def\picture #1 by #2 (#3){\leavevmode\vbox to #2{
     \hrule width #1 height 0pt depth 0pt
      \vfill
       \special{picture #3}}}

\def\illustration #1 by #2 (#3) scaled #4{\dimen1=#2
  \divide\dimen1 by 1000
  \multiply\dimen1 by #4
  \vtop to \dimen1{\dimen1=#1
  \divide\dimen1 by 1000
  \multiply\dimen1 by #4
  \hsize=\dimen1\vss
  \noindent\special{illustration #3 scaled #4}}}

\ifx\couleurs\oui

\fi

\anglais
\beginpackages
\endpackages
\optionparag=1
\def\paradouze{oui}
\vsize=255truemm
\voffset=-3truemm
\ifx\optionkeymacros\undefined\else \fi

\catcode`\Œ=\active\defŒ{{\aa}}       
\catcode`\º=\active\defº{\int}        
\catcode`\=\active\def{\c c}        
\catcode`\¶=\active\def¶{\partial}    
\catcode`\Ä=\active\defÄ{\oint}       
\catcode`\Æ=\active\defÆ{\triangle}   
\catcode`\Â=\active\defÂ{\neg}        
\catcode`\µ=\active\defµ{\mu}         
\catcode`\¿=\active\def¿{{\o}}        
\catcode`\¹=\active\def¹{\pi}         
\catcode`\Ï=\active\defÏ{{\oe}}       
\catcode`\§=\active\def§{{\ss}}       
\catcode`\ =\active\def {\dagger}     
\catcode`\Ã=\active\defÃ{\sqrt}       
\catcode`\·=\active\def·{\Sigma}      
\catcode`\Å=\active\defÅ{\approx}     
\catcode`\½=\active\def½{\Omega}      
\catcode`\£=\active\def£{{\it\$}}     
\catcode`\°=\active\def°{\infty}      
\catcode`\¤=\active\def¤{{\S}}        
\catcode`\¦=\active\def¦{{\P}}        
\catcode`\¥=\active\def¥{\bullet}     
\catcode`\»=\active\def»{\leavevmode\raise.585ex\hbox{\b a}}      
\catcode`\¼=\active\def¼{\leavevmode\raise.6ex\hbox{\b o}}        
\catcode`\­=\active\def­{\not=}       
\catcode`\²=\active\def²{\leq}        
\catcode`\³=\active\def³{\geq}        
\catcode`\Ö=\active\defÖ{\div}        
\catcode`\É=\active\defÉ{{\dots}}     
\catcode`\¾=\active\def¾{{\ae}}       
\catcode`\Ç=\active\defÇ{\og}         
\catcode`\Ò=\active\defÒ{``}          
\catcode`\Á=\active\defÁ{!`}          
\catcode`\¢=\active\def¢{\rlap/c}     
\catcode`\Ô=\active\defÔ{`}           
\catcode`\Õ=\active\defÕ{'}           


\catcode`\=\active\def{{\AA}}       
\catcode`\'=\active\def'{\c C}        
\catcode`\¯=\active\def¯{{\O}}        
\catcode`\¸=\active\def¸{\Pi}         
\catcode`\Î=\active\defÎ{{\OE}}       
\catcode`\®=\active\def®{{\AE}}       
\catcode`\×=\active\def×{\diamond}    
\catcode`\¡=\active\def¡{\accent'27}  
\catcode`\Ó=\active\defÓ{''}          
\catcode`\±=\active\def±{\pm}         
\catcode`\È=\active\defÈ{\fg}         
\catcode`\À=\active\defÀ{?`}          
\catcode`\Ð=\active\defÐ{--}          
\catcode`\Ñ=\active\defÑ{---}         


\catcode`\Š=\active\defŠ{\"a}        
\catcode`\'=\active\def'{\"e}        
\catcode`\•=\active\def•{\"{\i}}     
\catcode`\š=\active\defš{\"o}        
\catcode`\Ÿ=\active\defŸ{\"u}        
\catcode`\Ø=\active\defØ{\"y}        
\catcode`\å=\active\defå{\^A}        
\catcode`\€=\active\def€{\"A}        
\catcode`\…=\active\def…{\"O}        
\catcode`\†=\active\def†{\"U}        
\catcode`\‡=\active\def‡{\'a}        
\catcode`\Ž=\active\defŽ{\'e}        
\catcode`\'=\active\def'{\'{\i}}     
\catcode`\—=\active\def—{\'o}        
\catcode`\œ=\active\defœ{\'u}        
\catcode`\ƒ=\active\defƒ{\'E}        
\catcode`\æ=\active\defæ{\^E}        
\catcode`\é=\active\defé{\`E}        %
\catcode`\ˆ=\active\defˆ{\`a}        
\catcode`\=\active\def{\`e}        
\catcode`\"=\active\def"{\`{\i}}     
\catcode`\˜=\active\def˜{\`o}        
\catcode`\=\active\def{\`u}        
\catcode`\Ë=\active\defË{\`A}        
\catcode`\‹=\active\def‹{\~a}        
\catcode`\–=\active\def–{\~n}        
\catcode`\›=\active\def›{\~o}        
\catcode`\Ì=\active\defÌ{\~A}        
\catcode`\"=\active\def"{\~N}        
\catcode`\Í=\active\defÍ{\~O}        
\catcode`\‰=\active\def‰{\^a}        
\catcode`\=\active\def{\^e}        
\catcode`\"=\active\def"{\^{\i}}     
\catcode`\™=\active\def™{\^o}        
\catcode`\ž=\active\defž{\^u}        

\let\optionkeymacros\null

\def\EE{{\bb E}}
\def\gS{{\got S}}

\def\gc{{\got c}}

\font\tengp=cmmib10
\def\gth{\hbox{$\tengp\vartheta$}}
\def\sgth{\hbox{$\gp\sct\vartheta$}}
\font\ninegp=cmmib10 at 9pt
               \font\eightgp=cmmib8
              \font\sevengp=cmmib7
               
               \newfam\gpfam
               \textfont\gpfam=\tengp\scriptfont\gpfam=\sevengp
               \def\gp{\fam\gpfam\textfont1=\tengp\scriptfont1=\sevengp\tengp}%

\def\gphi{\hbox{$\gp\varphi$}}
\def\gth{\hbox{$\gp\vartheta$}}

\font\tenib=cmmib10
             \font\nineib=cmmib10 scaled 900
\font\sevenib=cmmib10 scaled 700
             \font\fiveib=cmmib10 scaled 500
\def\itg{\ib}
\def\gb{{\got b}}
\def\gz{{\got Z}}

\def\s{{\itg s}}

\def\som#1{\mathop{\sum^{#1}}}
\def\Vect{\mathop{\rm Vect}\nolimits}

\def\Vect{\mathop{\rm Vect}\nolimits}
\def\auteur{R. de la Bretche \& G. Tenenbaum}
 \def\titrart{On moments of the Erd\H os--Hooley Delta-function}
\hautspages{\auteur}{\titrart}
\dateheure
\titrecentre{\titrart}
\bigskip\medskip
\centerline {\auteur} 
\medskip
{\leftskip110mm
\obeylines
\it
To Roger Heath-Brown,
 on the occasion of 
 his seventy-fifth birthday
\par }

\medskip\bigskip
{\eightpoint\leftskip1cm\rightskip1cm
\noi{\bf Abstract.}  
For integer $n\geqslant 1$ and real $u$, let $\Delta(n,u):=|\{d:d\mid n,\,\e^u<d\leqslant \e^{u+1}\}|$. The Erd\H os--Hooley Delta-function is then defined by $\Delta(n):=\max_{u\in\r}\Delta(n,u).$ We provide new upper bounds for weighted real moments of this function.\PAR
\medskip
\noi
{\bf Keywords:}  Distribution of divisors, concentration of divisors, average order, Hooley's function, Erd\H os-Hooley's Delta-function, Waring's problem.\par
\smallskip 
\noi \bf 2020 Mathematics Subject Classification: \rm primary   11N37; secondary 11K65.\par }
\bigskip
\medskip

 \paraun{Introduction}
For integer $n\geqslant 1$ and real $u$, put 
  $$\Delta(n,u):=|\{d:d\mid n,\,\e^u<d\leqslant \e^{u+1}\}|,\quad\Delta(n):=\max_{u\in\r}\Delta(n,u).$$
 Introduced by Erd\H os \citer{Er73} (see also \citer{EN75}) and studied by Hooley \citer{Ho79}, the $\Delta$-function and specifically its (possibly weighted) mean-value proved very useful in several branches of number theory --- see, \eg,  \citer{Te86}, \citer{HT88}, and \citer{BT23}  for further references. If $\tau(n)$ denotes the total number of divisors of $n$, then $\Delta(n)/\tau(n)$ coincides with the concentration of the numbers $\log d$, $d|n$. \par 
 Asymptotic estimates for 
 $$S(x):=\sum_{n\leqslant x}\Delta(n)$$
 have a rather long history since Hooley's pioneer work \citer{Ho79}: see \citer{BT25} for references. Recent progress  regarding this mean-value is due to Koukoulopoulos and Tao \citer{KT23}, Ford, Koukoulopoulos  and Tao~\citer{FKT23}, and the authors \citer{BT25}. The current best  estimates are  $$x(\log_2x)^{1+\eta+o(1)}\ll S(x)\ll x(\log_2x)^{5/2},$$
 where $\eta\approx 0.353327$,  the lower and upper bounds being respectively proved   in \citer{FKT23} and \citer{BT25}. Here and in the sequel we let $\log_k$ denote the $k$-fold iterated logarithm.\par 
 In \citer{BT26}, we applied the latter upper bound to Waring type problems, showing that, for certain integer vectors $(2,\ell_1,\ldots,\ell_t)\in\N^{*(t+1)}$ satisfying $1/2+\sum_{1\leqslant j\leqslant t}1/\ell_j=1$,  the number of integers $n$ not exceeding $x$ and representable by the form $m_0^2+\sum_{1\leqslant j\leqslant t}m_j^{\ell_j}$ is $\gg x/(\log_2x)^{5/2}$.
 \par \medskip
 Given $z>0$, we  let $\M_z$ denote the class of those non-negative multiplicative functions $\varrho$ that are bounded on the set of prime powers  and satisfy, for suitable $c=c(\varrho)>0$,
$$\sum_{p\leqslant y}{ \varrho(p)\log p}=zy+O\Big({y\over (\log y)^c}\Big)\qquad (y\geqslant 2).\eqdef{condirho}$$
We then have
$$  \sum _{P^+(n)< x}{\mu(n)^2 \varrho(n)\over n}=\prod_{p<x}\Big(1+{\varrho(p)\over p}\Big)\asymp ({\log x } )^z\qquad (  x\geqslant 2). $$
Here and throughout, $\mu$ denotes the Mšbius function, the letter $p$ denotes a prime  number  and $\omega(n)$ denotes the number of distinct prime factors of an integer $n$. A standard  instance of an element of $\M_z$ is the function defined by $\varrho(n):=z^{\omega(n)}.$ 

 In this work, we investigate  upper bounds for the weighted moments
$$S_{t,\varrho}(x):=\sum_{n\leqslant x}\varrho(n)
\Delta(n)^t\qquad (t\geqslant 1,\,\varrho\in\M_z,\, x\geqslant 1).\eqdef{defStz}$$
Here and in the sequel, $t$ and $z$ are fixed parameters, $\varrho\in\M_z$, and implicit constants may depend upon $t$ and $\varrho$. 
 Estimates for $S_{t,\varrho}(x)$ are potentially useful for determining dominant values of $\varrho(n)\Delta(n)$ in certain averages. In \citer{Te90}, this was used with $t=1+\varepsilon$, for arbitrary small $\varepsilon$.\par 
Writing $ \beta= \beta(t,z):=2^tz-t$, $z_t:=t/(2^t-1)$, we have, for real $t\geqslant 1$, $$x(\log x)^{\max\{\beta,z\}-1}\ll S_{t,\varrho}(x)\ll x(\log x)^{\max\{\beta,z\}-1}\e^{\{\kappa+o(1)\}\sqrt{(\log_2 x)\log_3x}},\eqdef{HT}$$
with
$$\kappa:=\normalbaselineskip=15pt\cases{2t & if  $z\geqslant z_t$,\cr
2\sqrt{1-\log (1/z)/\log_22}& if  $z< z_t$.\cr}$$
The  upper bound is proved in \citer{HT86}\note{The proof of \citer{HT86} is written for $\varrho(n):=z^{\omega(n)}$ but it readily extends to $\varrho\in\M_z$.} while the
 lower bound follows from the inequality (see, \eg, \citeplus{HT88}{th.\thinspace60}) 
 $$\Delta(n)\geqslant \max\Big\{1,{\tau(n)\over 1+\log n}\Big\}\qquad (n\geqslant 1).\eqdef{triv}$$
  \medskip
Our approach to bounding $S_{t,\varrho}(x)$ is primarily based on an iterative procedure restricting the parameter $t$ to integer values.  While the induction step is based on the method developed in~\citer{KT23}, the initialization turns out to be  non trivial and necessitates new ideas---see  Propositions \ref{sEM2t/taut-ent} and~\ref{sSbfq} {\it infra}. The extension to real values of the exponent is derived secondarily by exploiting the flexibility offered by the parameter $z$.\par 
We obtain the following results. Here and in the sequel, we let $\delta=\delta(t,z):=\1_{\beta=z}=\1_{z=z_t}$.  
\Propt{thD}{Let  $t\geqslant 1$,  $z>0$, and  $\varrho\in\M_z$. \par 
\rm{(i)} For $z\geqslant z_t$ and all large  $x$, we have
$$ S_{t,\varrho}(x)
\ll x(\log x)^{\beta-1}(\log_2x)^{t+1+t\delta}.\eqdef{estsDt}$$
\par 
\rm{(ii)} For $0<z<z_t$ and all large $x$, we have
$$ S_{t,\varrho}(x)
\ll x(\log x)^{z-1}(\log_2x)^{2t+t/s},\eqdef{estsDts}$$
where $s>t$ is defined by $z_s=z$.
}\medskip
 \rems (i) As shown in \citer{BT25}, the upper bound  \eqref{estsDt} can be improved to $x(\log_2x)^{5/2}$ in the case $z=t=1$.
\par 
(ii) As will be clear from the proofs, all our estimates are locally uniform in $z$ and $t$ in the complement of the line $\beta=z$. In any admissible compact domain in $(z,t)$, they only depend on the constants involved in~\eqref{condirho} and on an upper bound for $\varrho$ on prime powers.
\smallskip
Inserting \eqref{estsDt} into the proof of \citeplus{NT98}{cor. 4} readily yields the following statement, where \hbox{\citeplus{BT25}{th.\thinspace1.1}} has been taken into account for the case $t=1$. 
\Propc{DQt}{Let $F\in\z[X]$ be irreducible with no fixed prime divisor. Then, for any $t\geqslant 1$ and~$\varepsilon>0$, we have
$$\sum_{x<n\leqslant x+y}\Delta\big(|F(n)|\big)^t\ll y(\log x)^{\beta(t,1)-1}(\log_2x)^{ t+1+\delta/2}$$
provided $x^\varepsilon\leqslant y\leqslant x$.}
\par \smallskip
For exponents at least equal to 2 and sufficiently large $z$, we obtain the true order of magnitude.
\vskip-3mm  
\Propt{thDelta3}{
 Let $t\geqslant 2$, $z>z^+_t:=t/(2^t-2^{t/2})=z_t(1+2^{-t/2})$, $\varrho\in\M_z$.
We have 
$$S_{t,\varrho}(x)\asymp x(\log x)^{\beta-1}\qquad (x\geqslant 2).\eqdef{Stzx}$$
}
We have $z_2^+=1$. This motivates the next statement.
  \Propt{thDelta2}{Let  $\varrho\in\M_1$. We have 
$$S_{2,\varrho}(x)=\sum_{n\leqslant x}\varrho(n)\Delta(n)^2\ll x(\log x)(\log_2x)\sqrt{\log_3x}\qquad (x\geqslant 16).\eqdef{S21x}$$}
 This last estimate improves  by a  factor essentially $\log_2x$ the current upper bound for the quadratic mean-value of the Delta-function: see 
\citeplus{BT25}{th. 1.2} where an argument from \citer{HT88} is made explicit.
\par 
 From \eqref{S21x} and \citeplus{NT98}{th. 1}, we derive the following corollary, which seems out of reach by standard techniques.
\Propc{corD2}{Let $\varepsilon\in]0,1[$. Uniformly for $N\geqslant 1$, $x\geqslant 16 $, $x^\varepsilon\leqslant y\leqslant x$, we have
$$\sum_{N<m,n\leqslant 2N}\bigg(\fl{x+y\over [m,n]}-\fl{x\over [m,n]}\bigg)\ll y(\log x)(\log_2x)\sqrt{\log_3x}.$$}
 \medskip
\paraun{Notation}
Put $$\eqalign{&E:=\{n\geqslant 1:\mu(n)^2=1\},\quad n_y:=\prod_{p|n,\,p<y}p\quad(n\geqslant 1,\,y\geqslant 1),\cr}$$
Given $t>0$ and $T\geqslant 2$,  put 
$$\leqalignno{& 
\gb=\gb_{t,z}:= 
{1\over 2^tz\log 2   -1}\qquad  (z\geqslant z_t),&\eqdef{b}\cr& f_T(y):=\gc\min\bigg(\log (T\log 3y),{(\log_23y-\gb_{t,z}\log T)^2\over \log T}\bigg)\quad(y\geqslant 1),&\eqdef{fT}\cr}$$
where $\gc$ is a strictly positive constant to be chosen sufficiently small in the sequel.
\par 
Define
$$\eqalign{ E_T^{t,z}&:=\big\{n \in E :\tau(n_y)\leqslant \,\e^{-f_T(y)}T\log 3y \ (y\geqslant 1)\big\},\cr
E_T &:=\big\{n \in E :\tau(n_y)\leqslant T\log 3y \ (y\geqslant 1)\big\}.}\eqdef{defETt}$$
 It is useful to bear in mind that $n\in E_T^{t,z}$ implies that $d\in E_T^{t,z}$ for any divisor $d$ of $n$.
\par 
Let $P^+(n)$---resp. $P^-(n)$---denote the largest---resp. the smallest---prime factor of $n>1$ with the convention that $P^+(1):=1$, $P^-(1)=\infty$.
Given $x\geqslant 2$, we consider the probability $\PP_{x,\varrho}$  defined on $E$~by $$\PP_{x,\varrho} (\{n\}):={\varrho(n)\over n}\prod_{2\leqslant  p<x}\Big(1+{\varrho(p)\over p}\Big)^{-1}\asymp{\varrho(n)\over n  ( \log x)^z } \qquad \big(P^+(n)<x\big),$$
and let $\EE_{x,\varrho}$ denote the corresponding expectation. 
\par 
Throughout this article, the symbol $\som{x}$ indicates a summation over squarefree integers whose prime factors  are restricted to the interval $[2,x[$.
For the sake of future reference, we note that, if~$\varrho\in\M_z$, then $$\som{x}_{n\geqslant 1}{\varrho(n)\over n}\asymp (\log x)^z\qquad (x\geqslant 2).\eqdef{condirho2}$$
\par 
\medskip\goodbreak
\paraun{Basic lemmas}
 \vskip-5mm
 \Propl{2.1}{Let $t\geqslant 1$, $z>0$  and  $\varrho\in\M_z$. We have  $$ S_{t,\varrho}(x) \asymp x (\log 3x)^{z-1}\,\EE_{x,\varrho}(\Delta^t)\qquad (x\geqslant 1). \eqdef{Marith-log}$$}
 \nid This is a trivial extension of \citeplus{HT88}{th.\thinspace61}. We omit the details. 
\qed
\smallskip
As in \citer{KT23}, the first step to the proof of \ref{thD} consists in defining a set of integers with useful multiplicative constraints.
 
 \par 

The next lemma is analogous to \citeplus{KT23}{prop. 4.1}. The essential feature consists in bounding $\PP_{x,\varrho}( E\sset E_T^{t,z})$ by a multiple of $1/T^{ t}$. The fact that the resulting estimate is trivial for small~$T$ will have no consequence.
\par  
Recall notation $z_t:=t/(2^t-1)$, so that $z=z_t$ if, and only if, $\beta=z$.
\Propl{PE-ET}{Let $t\geqslant 1 $, $z\geqslant z_t$. We have 
$$\PP_{x,\varrho}\big( E\sset E_T \big) \leqslant  \PP_{x,\varrho}\big( E\sset E_T^{t,z}\big) \ll {(\log x)^{\beta-z}\over T^ t}\qquad (x\geqslant 3,\,T\geqslant 3).\eqdef{EssetET}$$}

\nid 
We may assume $ T^t\geqslant (\log 3x)^{\beta-z}$ since the stated estimate is otherwise trivial. Put $$\kappa_{y,z}:={1\over \log 2}\Big(1+{\log T-f_T(y)\over \log_23y}\Big),$$
and note right away that, since $f_T(y)\leqslant \gc\log T+\gc\log_23y$, we have
$$\kappa_{y,z}\geqslant {t+(1-\gc)(\beta-z)-\gc t\over t\log 2}={(1-\gc)(2^t-1)z\over t\log 2}={(2^t-1)z\over t\log 2}-O(\gc)>z,\eqdef{minkappay}$$
provided $1\leqslant y\leqslant x$ and $\gc$ is chosen sufficiently small.
\par \goodbreak
We have $$E_{T,x}^{t,z}:=E_T^{t,z}\cap\{n\geqslant 1:P^+(n)<x\}=\{n\in E:P^+(n)<x,\,\sup_{y\geqslant 1}\{\omega(n_y)/(\kappa_{y,z}\log_23y)\}\leqslant  1\}.$$ 
Indeed, $\log T+\log_2(3y)-f_T(y)$ is increasing in $y$ for
sufficiently small $\gc$. As a function of $\log_2(3y)$, it is also
Lipschitz with a constant independent of $T$. To check both assertions,
note that on the quadratic branch defining $f_T$, the inequality
$$\{\log_2(3y)-\gb\log T\}^2
 \leqslant\log T\{\log T+\log_2(3y)\}$$
bounds $\log_2(3y)/\log T$ in terms of $\gb$.\par 
If $n\in E\sset E_{T,x}^{t,z}$ and $P^+(n)<x$, then  $\omega(n_y)\geqslant \kappa_{y,z} \log_23y-c$ for a suitable constant $c$   independent of $T$, $y$, $x$ and all $y$ in a suitable interval $[y_0,y_0^2]$  with $y_0\leqslant x$. Therefore,
$$\eqalign{\1_{E\sset E_{T,x}^{t,z}}(n)&\ll \int_{y_0}^{y_0^2}\sum_{k\geqslant \kappa_{y,z} \log_23y-c}\1_{\omega(n_y)=k}{\dd y\over y\log 3y},\cr
\PP_{x,\varrho}(E\sset E_{T}^{t,z})&\ll{1\over (\log 3x)^{z}}\int_1^{x^2}
\sum_{k\geqslant \kappa_y\log_23y-c}{1\over k!}\Big(\sum_{p\leqslant y}{\varrho(p)\over p}\Big)^k\som{y,x}_{n\geqslant 1}{\varrho(n)\over n}{\dd y\over y\log 3y}\cr
&\ll \int_1^{x^2}{\e^{c\kappa_{y,z}}\dd y\over y(\log 3y)^{1+zQ(\kappa_{y,z} /z)}\sqrt{\log_23y}},\cr}$$
where $Q(v):=v\log v-v+1$ $(v>0)$. 
Here the sum over $k$ has been estimated by a standard bound for partial sums of the exponential series such as \citeplus{No76}{lemma 4.7}, taking \eqref{minkappay} into account and observing that the perturbations of the Poisson parameter and the threshold contribute $\ll \kappa_{y,z}$ to the exponent.   Now observe that $$Q(2^t)=2^t\log 2^t-2^t+1,\quad Q'(2^t)=\log 2^t,\quad 
 Q''(v)\gg1/v$$ for some $v=\vartheta2^t+(1-\vartheta)\kappa_{y,z} /z$ with $0\leqslant \vartheta\leqslant 1$. Therefore, there is a constant $\gc_0=\gc_0(z,t)>0$ such that
$$\eqalign{1+zQ(\kappa_{y,z} /z)&\geqslant 1+zQ(2^t)+ (\kappa_{y,z} -2^tz)Q'(2^t)+  {\gc_0\big(\kappa_{y,z} -2^t z\big)^2\over \kappa_{y,z}+2^tz}
\cr
&=1+z-2^tz+\kappa_{y,z}  t\log 2+ {\gc_0  (\kappa_{y,z} -2^t z)^2\over \kappa_{y,z}+2^tz}\cr
&=1+z-\beta+{t(\log T-f_T(y))\over \log_23y}+ {\gc_0  (\kappa_{y,z} -2^t z)^2\over \kappa_{y,z}+2^tz},\cr}$$
whence
$$ \e^{-c\kappa_{y,z}}(\log 3y)^{1+zQ(\kappa_{y,z} /z)}\geqslant T^t(\log 3y)^{1+z-\beta}\e^{ -c\kappa_{y,z}-  t f_T(y)+ \gc_0(\kappa_{y,z} -2^tz)^2(\log_23y)/(\kappa_{y,z}+2^tz)}.$$
since, for sufficiently large $y$, the term $-c\kappa_{y,z}$ is absorbed into the implicit constant when $\kappa_{y,z}$ is bounded and into half the Taylor remainder otherwise.
For suitable $\gc_1>0$, the argument of the latter exponential is $\geqslant \gc_1\log T$ in the range $\log 3y\leqslant T^{\gb/2}$, while it is $\geqslant \gc_1\log_23y$ for $\log 3y>T^{2\gb}$.  Therefore, after excluding bounded values of $y$,  which contribute $\ll1/T^t$ since the corresponding events are empty for sufficiently large $T$, we see that the contribution of the union of these ranges to the integral is 
$$\ll{(\log x)^{\beta-z}\over T^{t}}\Big({\log T\over T^{\gc_1t}}+{1\over (\log T)^{2\gc_1\gb t}}\Big).$$
\par 
In the remaining middle range, $\kappa_{y,z}$ stays in a fixed compact subinterval of $]z,\infty[$. However, writing $D:=\gb\log T-\log_23y$, we have
$$(\kappa_{y,z}-2^tz)\log_23y
 ={D\over\gb\log 2}-{f_T(y)\over\log 2}\cdot$$
If $D\geqslant0$, then $D\leqslant\gb\log T$ and
$$f_T(y)\leqslant{\gc D^2\over\log T}\leqslant\gc\gb D.$$
Thus, choosing $\gc$ sufficiently small, we obtain
$$\big|(\kappa_{y,z}-2^tz)\log_23y\big|
 \geqslant {|D|\over2\gb\log 2}\cdot$$
The same inequality holds when $D<0$, since $f_T(y)\geqslant0$.
Moreover,
$$\eqalign{
(\kappa_{y,z}+2^tz)\log_23y
 &= {\log T-f_T(y)+\log_23y\over\log 2}
       +2^tz\log_23y\ll\log T+\log_23y.\cr}$$
Consequently,
$${(\kappa_{y,z}-2^tz)^2\log_2(3y)\over\kappa_{y,z}+2^tz}
 \gg {D^2\over\log T+\log_23y}\cdot$$
On the other hand, the definition of $f_T$ implies
$$f_T(y)\ll{\gc D^2\over\log T+\log_23y}\cdot$$
Indeed, if $\log_23y\leqslant2\gb\log T$, this follows from
$f_T(y)\leqslant\gc D^2/\log T$. Otherwise,
$|D|\geqslant{1\over2}\log_23y$, and we use
$f_T(y)\leqslant\gc\{\log T+\log_23y\}$.
Taking $\gc$ sufficiently small therefore gives
$${\gc_0(\kappa_{y,z}-2^tz)^2\log_23y\over\kappa_{y,z}+2^tz}
 \geqslant(t+1)f_T(y).\eqdef{mincarre}$$
The remaining integral is therefore
$$\ll{(\log3x)^{\beta-z}\over T^t}
 \int_{\di{y\geqslant2}{T^{\gb/2}\leqslant\log 3y\leqslant T^{2\gb}}}
 {\e^{-f_T(y)}\d y\over y\log3y\sqrt{\log_23y}}
 \ll{(\log3x)^{\beta-z}\over T^t}\cdot$$
Indeed, in this range $\log_23y\asymp\log T$ and
$f_T(y)\gg\{\log_2(3y)-\gb\log T\}^2/\log T$, so the last
integral is bounded by a Gaussian integral. This completes the proof.
\qed
\smallskip
\rem The factor $\e^{-f_T(y)}$ appearing  in the definition \eqref{defETt} of the set $E_T^{t,z}$ will only be useful in the case $\beta=z$, \ie $z=z_t$. However, this will turn out to be crucial in the proof of part (ii) of \ref{thD}.  
\par 
 \par 
 \bigskip
 \paraun{A functional inequality} Denote by $\P$ the set of all prime numbers.
Consider a non-negative arithmetic function~$f$ supported on the set of squarefree integers. Assume furthermore that $f(1)=1$ and that, for a suitable function $w:\N^*\times \P\to\r^+$ vanishing if the first variable is not squarefree, we have 
$$f(mp)\leqslant f (m) +w (m,p) \qquad \big(  m\geqslant 1,\,p\nmid m \big).\eqdef{relf}$$
\par 
The following lemma furnishes a useful bound for
$${\scal F}(x):= \som x_{n\geqslant 1}{\varrho(n)f(n)\over n}\cdot$$
It has been established in  \citer{KT23} for $\varrho=\1$, and hence $z=1$---see also \citer{BT25}.
 \Propl{lemmerhof}{Let  $z>0$ and $\varrho\in\M_z$.  For $f,\,w$ satisfying \eqref{relf}, we have $${{\scal F}(x)\over (\log x)^z}\ll 1+  \int_1^{x^2}\som y_{m\geqslant 1}{\varrho(m)\over m }\sum_{\sqrt{y}\leqslant p< y}{\varrho(p)w(m,p)\over p}{\d y\over y( \log 3y )^{1+z}}\cdot\eqdef{inegfonct}$$
}

For $n\geqslant 1,\,q\geqslant 1,$ define
$$\leqalignno{M_q(n)&:=\int_\r\Delta(n,u)^q\d u,&\eqdef{defMq}\cr
\mu_q(n)&:=M_q(n)/\tau(n).&\eqdef{defmuq}\cr}$$
We shall show that \eqref{inegfonct} holds  for $f=\1_S\mu_q^r$ with $ 1\leqslant r\leqslant t $ and $q\geqslant 3$,	where  $S$ is some set of integers such that $(mp\in S)\Rightarrow (m \in S)$.  
\medskip\goodbreak
 \nid  Writing $n=mp$ with  $P^+(m)<p$, we get
$${\scal F} (x)\leqslant 1+\sum_{p<x} \som p_{m\geqslant 1}{\varrho(p)\varrho(m)f(pm)\over pm}\cdot$$
 Applying \eqref{relf}  when $\mu(pm)^2=1$   yields  
$${\scal F} (x) \leqslant 1+\sum_{p<x}\som p_{m\geqslant 1} {\varrho(p)\varrho(m)(f(m)+w (m,p))\over pm }= \sum_{p<x} {\varrho(p){\scal F} (p)\over p }+V(x),$$
with
$$V(x):=1+\sum_{p<x}{\varrho(p) \over p  } \som p_{m }{ \varrho(m)w(m,p)\over  m }\cdot$$
Iterating this as in \citer{KT23} and taking \eqref{condirho2} into account, we get
$${\scal F} (x)\leqslant V(x)+\sum_{2\leqslant P(n)<x} {\mu^2(n)\varrho(n)V(P^-(n))\over n}\ll 
V(x)+\sum_{p<x} {\varrho(p)V(p)\over p}\bigg({\log x\over \log p}\bigg)^z,$$
whence
$${{\scal F}(x)\over (\log x)^z}\ll {V(x)\over (\log x)^z}+\sum_{p<x}{\varrho(p)V(p)\over p(\log p)^z}\cdot\eqdef{eqfonctT}$$
Indeed, at the first step we obtain
$$\sum_{p<x}{\varrho(p){\scal F}(p)\over p}\leqslant \sum_{p_1<x}{\varrho(p_1)\over p_1}\bigg\{\sum_{p_2<p_1}{\varrho(p_2){\scal F}(p_2)\over p_2}+V (p_1)\bigg\},$$
 then
$$\eqalign{\sum_{p<x}{\varrho(p){\scal F} (p)\over p}&\leqslant \sum_{p_1<x}{1\over p_1}\bigg(\varrho(p_1)V(p_1)+\sum_{p_2<p_1}{\varrho(p_1p_2)\over p_2}\Big\{\sum_{p_3<p_2}{\varrho(p_3){\scal F} (p_3)\over p_3}+V(p_2)\Big\}\bigg)\cr
&\leqslant \sum_{p_1<x}{\varrho(p_1)V(p_1)\over p_1}+\sum_{p_2<p_1<x}{\varrho(p_1p_2)V(p_2)\over p_1p_2}+\sum_{p_3<p_2<p_1<x}{\varrho(p_1p_2p_3){\scal F}(p_3)\over p_1p_2p_3},\cr}$$
and so on.
\par Since $n_y=n_p$ when  {$P^+(n)<p<y\leqslant p^2$}, we may write,  still following \citer{KT23},
$$\eqalign{V(x)&\ll 1+\sum_{p<x}\som p_{m }{\varrho(p)\varrho(m)w (m,p)\over pm} \int_p^{p^2}{\d y\over y\log 3y}\cr
&\ll 1+ \int_1^{x^2}\som y_{m }{\varrho(m)\over m }\sum_{\sqrt{y}\leqslant p< y}{\varrho(p)w(m,p)\over p}{\d y\over y \log 3y }\cdot
\cr}$$  
Inverting summation and integration in \eqref{eqfonctT} and appealing to the uniform bound
$$\sum_{\sqrt{y}\leqslant p<x}{\varrho(p)\over p(\log p)^z}\ll  {1\over (\log y)^z},$$ 
stemming from \eqref{condirho}, we obtain  \eqref{inegfonct} as stated.
\qed
\medskip
\bigskip
\paraun{Moments of quadratic moments}\drefun{E*T}
 \par 
 Recall the notation $\delta =\delta(t,z):=\1_{\beta=z}=\1_{z=z_t}$. 
 \Propp{EM2t/taut-ent}{Let $z>0$ and  $\varrho\in\M_z$.
 \par 
 {\rm(i)} If  $t\in\r^{+*}$, $ z\geqslant z_t $, we have  
$$ \EE_{x,\varrho}\big(  \mu_2^t\big) \gg (\log x)^{\beta-z} 
\qquad (x\geqslant 3 ).\eqdef{EEM2/taurgeq}$$
\par
 {\rm(ii)} If $t\in\N^*$, $ z\geqslant z_t $, we have  $$\leqalignno{ \EE_{x,\varrho}\big(  \mu_2^t\big) &\asympÊ(\log x)^{\beta-z}(\log_2x)^\delta
\qquad (x\geqslant 3 ),&\eqdef{EEM2/taur}\cr
\EE_{x,\varrho}\big(\1_{E_T^{t,z}}  \mu_2^{t+1}\big) &\ll T (\log x)^{\beta-z}
\qquad (T\geqslant 3,\,x\geqslant 3 ).&\eqdef{EEM2/taut+1}}$$
If moreover  $z>z_t$, then 
$$\EE_{x,\varrho}\big(\1_{E_T}  \mu_2^{t+1}\big) \ll T (\log x)^{\beta-z}
\qquad (T\geqslant 3,\, x\geqslant 3 ).\eqdef{EEM2/taut+1*}
$$
 {\rm(iii)} If $t\geqslant 1$, $z>z_t$, we have
$$\EE_{x,\varrho}\big(\mu_2^t\big)\asymp(\log x)^{\beta-z}.\eqdef{mu2tgen}$$}
  \nid
  (i) Put $\tau(n,\vartheta):=\sum_{d|n}d^{i\vartheta}$ $(n\geqslant 1,\vartheta\in\r)$. By Parseval's formula (see \citeplus{HT88}{(6.23)})
we have $$\mu_2(n)\asymp \int_0^{1}{|\tau(n,\vartheta)|^2\over \tau(n)}\d\vartheta.\eqdef{Pars}$$
Estimate \eqref{EEM2/taurgeq}  follows from the observation that, for $x\geqslant 3$, $\mu(n)^2=1$, $P^+(n)<x$, we have
$${|\tau(n,\vartheta)|^2\over \tau(n)}=\prod_{p|n}\{1+\cos(\vartheta\log p)\}\geqslant \tau(n)\prod_{p|n}\Big(1-{(\log p)^2\over 4(\log x)^2}\Big)\qquad \big(0\leqslant \vartheta\leqslant 1/\log x\big).$$
Indeed, we  infer from this and \eqref{Pars} that
$$\eqalign{\EE_{x,\varrho}\big(\mu_2^t\big)&\gg {1\over (\log x)^{z+t}} \som x_{n\in E}{\varrho(n)\tau(n)^t\over n}\prod_{p|n}\Big(1-{(\log p)^2\over 4(\log x)^2}\Big)^t\cr
&
\gg {1\over (\log x)^{z+t}} \exp\bigg\{ \sum_{p\leqslant x} {\varrho(p)2^t\over p}\Big(1-{(\log p)^2\over 4(\log x)^2}\Big)^t\bigg\}\asymp{(\log x)^{\beta-z}},
}$$
where we used \eqref{condirho} to estimate the last sum over $p$.
\par
(ii) First apply \eqref{Pars} to get
$$\eqalign{\EE_{x,\varrho}\big(\mu_2^t\big)&\asymp {1\over (\log x)^z}\int_{\sgth \in [0,1]^t}Ê\som{x}_{n\geqslant 1}{\varrho(n)\over n} \prod_{1\leqslant j\leqslant
t}{|\tau(n;\vartheta_j)|^2 \over \tau(n)} \dd\gth \cr
&\asymp{1\over (\log x)^z}\int_{\sgth \in [0,1]^t}\exp\Bigg\{\sum_{p\leqslant x}{\varrho(p)\over
p }\prod_{1\leqslant j\leqslant
t}{|\tau(p;\vartheta_j)|^2\over 2}\Bigg\}\dd\gth 
 .\cr}\eqdef{EMr/taur}$$
\par 
At this stage some notation is necessary. We use those of \citer{BT12}. Given an integer
$s\geqslant 1$, we designate by
$\{e_j:1\leqslant j\leqslant s\}$ the canonical basis of $\r^s$ and by  $\W_s$  the set of linear forms $w\in\L^*(\r^s):=\L(\r^s,\r)$
such that $w(e_j)\in\{-1,0,1\}$ for all $j$.  Note that $0_{\L^*(\r^s)}\in \W_s$. We define the length of $w\in\W_s$ by  $$|w|:=\sum_{1\leqslant j\leqslant s} |w(e_j)|.\eqdef{defL1}$$

Given a prime number $p$ and $\gth\in \r^t$, we have
$$\eqalign{\prod_{1\leqslant j\leqslant
t}{|\tau(p;\vartheta_j)|^2\over \tau(p)}&=\prod_{1\leqslant
j\leqslant
t}\big(1+\dm(p^{i\theta_j}+p^{-i\theta_j})\big)
=\sum_{w\in\W_t} 
{\cos\big(w(\gth)\log
p\big)\over 2^{|w|} }.\cr}$$
In view of the estimate
$$
\sum_{p\leqslant y}{\varrho(p)\cos (\psi \log p)\over p} =z\log \Big({\log y\over 1+\psi\log y}\Big)+O(1)\quad(y\geqslant 2,\,0\leqslant \psi\leqslant 1),\eqdef{estimatesump}
$$
that follows from \eqref{condirho} by partial summation,
we derive that
$$\eqalign{\sum_{p\leqslant x}{\varrho(p)\over p}\prod_{1\leqslant j\leqslant
t}{|\tau(p;\vartheta_j)|^2\over 2}
&=\sum_{w\in\W_t} {1\over 2^{|w|}}\sum_{p\leqslant
x} {\varrho(p)\cos\big(w(\gth)\log
p\big)\over p}\cr&=
 \sum_{w\in\W_t}{z\over 2^{|w|}}\log
\Big(1+{\log x\over 1+|w(\gth)|\log x}\Big)
+O(1).\cr}$$

Defining
$$\eqalign{I_{t,z}(X) :=\int_{[0,X]^t}
&\prod_{w\in \W_t} \Big(1+{X\over
1+|w(\gth)|
}\Big)^{z/2^{|w|}} {\dd\gth},
\cr}\eqdef{intX}$$
we get, after homothetical change of variables, 
$$\EE_{x,\varrho}\big( \mu_2^t\big) \asympÊ{I_{t,z}(\log x)\over (\log x)^{t+z}}\cdot $$
It hence remains to evaluate $I_{t,z}(X)$, bearing in mind that
$$\sum_{w\in \W_t}{1\over 2^{|w|}}=\sum_{0\leqslant j\leqslant t}{1\over 2^j}\sum_{\di{{w\in \W_t}}{|w|=j}}1 = \sum_{0\leqslant j\leqslant t} {t\choose j}=2^t.\eqdef{inegWt}$$
\par 
We claim that
$$I_{t,z}(X) \asymp  X^{2^tz}\int_{[0,X]^t}
 \prod_{w\in \W_t} \Big({1\over
1+|w(\gth)|
}\Big)^{z/2^{|w|}} {\dd\gth} \asymp  X^{2^tz}(\log X)^{\delta},\eqdef{majIt}$$
 which plainly implies \eqref{EEM2/taur}.
\par 
Let us first establish  the lower bound included in \eqref{majIt}. Since we plainly have  $w(\gth)\ll \vartheta_t$  if $$0<\dm \vartheta_t/4^{t-j}\leqslant \theta_j\leqslant  \vartheta_t/4^{t-j}\quad(1\leqslant j\leqslant t)$$ and since the subsum  
of \eqref{inegWt} corresponding to non-zero  forms equals $2^t-1$, the integral is  $$\gg  \int_{1}^X{\dd\vartheta_t\over \vartheta_t^{(2^t-1)z-(t-1)}}
\cdot$$ 
 This yields the lower bound included in \eqref{majIt}.
\par 
 We  now prove the  upper bound included in \eqref{majIt}. We may assume $z_t\leqslant z<z_{t-1}$ since the  integrand is a non-increasing function of  $z$. 
The first step consists in associating to each $\gth\in ]0,+\infty[^t$ a basis~$\B_{\sgth}$ of $\L(\r^t,\r)\cap\W_t$ constructed in the following way : $w_1$ is the non-zero linear form of~$\W_t$ minimising  $|w(\gth)|$,  $w_2$   minimises $|w(\gth)|$ on $\W_t\sset\Vect (w_1)$, and for each index $k\in[1,t]$, the form~$w_k$   minimises $|w(\gth)|$ on $\W_t\sset\Vect (w_1,w_2,\ldots, w_{k-1})$.

The set $B:=\{\B\in\W_t^t:\exists\gth\in]0,\infty[^t:\B=\B_{\sgth}\}$ is finite.
Let us assume temporarily that, for any basis $\B=(w_1,w_2,\ldots, w_t)$ in $B$, we have  $$c_k=c_k(\B) :=\sum_{\di{w\in \W_t\cap\Vect (w_1,w_2,\ldots, w_{k}) }{w\neq 0}}{1\over 2^{|w|}}\leqslant 2^k-1\qquad (1\leqslant k\leqslant t).\eqdef{inegck}$$
This is shown in \ref{lemmeck} below,  which also yields $c_t=2^t-1$. \par 
Given a basis $\B$ in $B$, consider the domain $$\D(\B):=\{\gth\in[0,X]^t:\B_{\sgth}=\B\}.\eqdef{defDB}$$ 
 We have
$$\prod_{w\in \W_t} \Big({1\over
1+|w(\gth)|
}\Big)^{z/2^{|w|}}\ll \prod_{1\leqslant j\leqslant t} \Big({1\over
1+|w_j(\gth)|
}\Big)^{z(c_j-c_{j-1})}\qquad \big(\gth\in\D(\B)\big),
$$
where $c_j:=c_j(\B)$ $(1\leqslant j\leqslant t)$ and, by convention, $c_0=0$.
Performing  the change of variables  $\gphi=(w_1(\gth), \ldots, w_t(\gth))$, and noting that the associated Jacobian is constant, we get  
$$\int_{\D(\B)}\prod_{w\in \W_t} \Big({1\over
1+|w(\gth)|
}\Big)^{z/2^{|w|}} {\dd\gth}\ll\int_{|\varphi_1|\leqslant \ldots\leqslant |\varphi_t|\leqslant t X}\prod_{1\leqslant j\leqslant t}  \Big({1\over
1+|\varphi_j |
}\Big)^{z(c_j-c_{j-1})}\d\gphi.\eqdef{intDB}$$ 
Integrate successively  according to variables $\varphi_1$, $\varphi_2$, etc. 
If $z_t\leqslant z<z_{t-1}$, we use $zc_j<j$ for $j<t$. The last integral is
$$\int_0^{tX}{\dd \varphi_t\over (1+\varphi_t)^{\beta-z+1}}\ll(\log 2X)^\delta.$$
If $z\geqslant z_{t-1}$, the integrand does not exceed that of a $z'\in]z_t,z_{t-1}[$. The case $t=1$ corresponds to the integral $\int_0^X\d\varphi/(1+\varphi)^z$. This completes the proof of \eqref{majIt}.

\par 
It remains to establish \eqref{EEM2/taut+1} and \eqref{EEM2/taut+1*}.
We give the details both at the critical value $z=z_t$ and for $z>z_t$.
Put $d=t+1$. Throughout this part of the proof $t$ is a
positive integer, $x\geqslant3$, and $T\geqslant3$. We use the notation
$S=E_T^{t,z}$ and $F=f_T$ when $z=z_t$, and $S=E_T$ and $F=0$ when
$z>z_t$. Thus, in either case,
$$\tau(n_y)\leqslant T(\log 3y)\e^{-F(y)}\qquad(n\in S,\ y\geqslant1).$$
By \eqref{Pars}, we have
$$\EE_{x,\varrho}\bigg(\1_S\mu_2^{t+1}\bigg)
\ll {1\over (\log x)^z}\int_{[0,1]^d}
\som x_{n\in S}{\varrho(n)\over n}
\prod_{1\leqslant j\leqslant d}
 {|\tau(n;\vartheta_j)|^2\over\tau(n)}\,\dd\gth.
\eqdef{EM2t/taut+}$$
\par 
As above, partition the integration domain according to the basis
$\B=(w_1,\ldots,w_d)$ selected by successive minima in $\W_d$.
For each such basis put
$$\eqalign{
\D(\B)&=\{\gth\in[0,1]^d:\B_{\sgth}=\B\},\quad
\varphi_j=w_j(\gth),\qquad s_j=|\varphi_j| \quad(1\leqslant j\leqslant d),\cr
\qquad c_0&=0, \quad c_j=\sum_{\di{w\in\W_d\cap\Vect(w_1,\ldots,w_j)}{w\ne0}}
             2^{-|w|}\quad(1\leqslant j\leqslant d)\quad(1\leqslant j\leqslant d),.
}$$
Then $0\leqslant s_1\leqslant\cdots\leqslant s_d\leqslant d$,
$c_d=2^d-1$, and \ref{lemmeck} gives $c_j\leqslant2^j-1$.
The finitely many changes of variables have constant, non-zero
Jacobians; splitting into sign domains consequently costs only a
constant depending on $t$.
To include small and large frequencies without a further change of
variables, define
$$R(s)=\max\{1/\log x,\min(1,s)\},\qquad
\lambda_j={1\over R(s_j)},\qquad y_j=\e^{\lambda_j}.
$$
Thus $1\leqslant\lambda_d\leqslant\cdots\leqslant\lambda_1
\leqslant \log x$, and $\e\leqslant y_j\leqslant x$.
In particular, the cutoffs are now defined on the whole integration
domain. Write
$$U_p(\gth)=\prod_{1\leqslant j\leqslant d}
 {|\tau(p;\vartheta_j)|^2\over2}
=\sum_{w\in\W_d}{\cos\big(w(\gth)\log p\big)\over 2^{|w|}},
\qquad b_p={\varrho(p)U_p(\gth)\over
 p\sqrt{\tau(p_{y_t})\tau(p_{y_d})}}.$$
The numbers $b_p$ are non-negative and $b_p\ll1/p$, uniformly in
$\gth$ and $x$.
\par 
We first justify the saving obtained by prescribing the number of
prime factors below $y_d$. Give each prime $p<x$ an independent
Bernoulli variable with success probability $b_p/(1+b_p)$, and let
$K$ be the sum of these variables over $p<y_d$. Then
$$\sup_{r\in\z}\PP(K=r)\ll{1\over\sqrt{\log(2\lambda_d)}}\cdot$$
Indeed, when $\lambda_d$ is large, all coordinates $\vartheta_j$ are
$\ll_t1/\lambda_d$. For a sufficiently small fixed $\varepsilon>0$,
we have $U_p(\gth)\gg_t1$ for $p\leqslant\exp(\varepsilon\lambda_d)$.
It follows from \eqref{condirho} that
$$V:=\sum_{p<y_d}{b_p\over(1+b_p)^2}
\gg\log\lambda_d-O(1).$$
By Fourier inversion, every atom is
$$\ll {1\over2\pi}\int_{-\pi}^{\pi}
 \e^{-2V\sin^2(v/2)}\d v\ll{1\over \sqrt{1+V}}\cdot$$
The assertion for bounded $\lambda_d$ follows from the trivial bound
for a probability. All constants are independent of the prescribed
integer $r$.
\par 
For $j=t,d$ put $A_j=T(\log3y_j)\e^{-F(y_j)}$ and split the sum over
$n\in S$ into the classes
$$2^{-h-1}<{\tau(n_{y_t})\over A_t}\leqslant2^{-h},\qquad
  2^{-\ell-1}<{\tau(n_{y_d})\over A_d}\leqslant2^{-\ell}
  \qquad(h,\ell\geqslant0).$$
An empty class contributes nothing. In a non-empty class, the second
condition fixes $\omega(n_{y_d})$ to a single integer, since
$\tau(n_{y_d})=2^{\omega(n_{y_d})}$. Multiplication by
$$ {\sqrt{A_tA_d}\over
  2^{(h+\ell)/2}\sqrt{\tau(n_{y_t})\tau(n_{y_d})}}$$
changes each summand by a factor between $1$ and $2$.
Dropping all restrictions except the prescribed value of
$\omega(n_{y_d})$, applying the preceding atom bound, and summing
the convergent geometric series in $h,\ell$, we obtain
$$\eqalign{
\som x_{n\in S}{\varrho(n)\over n}
 \prod_{1\leqslant j\leqslant d}
 {|\tau(n;\vartheta_j)|^2\over\tau(n)}
&\ll {T\sqrt{\lambda_t\lambda_d}\,\e^{-F(y_d)/2}
       \over\sqrt{\log(2\lambda_d)}}
       \e^{\sum_{p<x}b_p}.
}$$
Here we used $F\geqslant0$ and $\log3y_j\asymp\lambda_j$.
\par  
We next estimate the attenuated Euler product. Set
$$\alpha:=1-{1/\sqrt2},\qquad a_t=\alpha(1+c_t).$$
The attenuation factors are $1/2$ for $p<y_d$, $1/\sqrt2$ for
$y_d\leqslant p<y_t$, and $1$ for $p\geqslant y_t$. Hence
$$\eqalign{
\sum_{p<x}b_p
={}&\sum_{p<x}{\varrho(p)U_p(\gth)\over p}
-{1\over2}\sum_{p<y_d}{\varrho(p)U_p(\gth)\over p}-\alpha\sum_{y_d\leqslant p<y_t}
                  {\varrho(p)U_p(\gth)\over p}.
}$$
Formula \eqref{estimatesump}, with bounded frequencies and the
regularized cutoffs above, provides
$$\eqalign{
\sum_{p<y_d}{\varrho(p)U_p(\gth)\over p}
 &=z2^d\log\lambda_d+O(1),\cr
\sum_{y_d\leqslant p<y_t}{\varrho(p)U_p(\gth)\over p}
 &=z(1+c_t)\log(\lambda_t/\lambda_d)+O(1).
}$$
For the second estimate, the forms in $\Vect(w_1,\ldots,w_t)$ have
frequency $O_t(s_t)$ and total mass $1+c_t$, including the zero
form. Every other form has frequency at least $s_d$, by the
successive-minimum construction, and contributes $O(1)$ on this
interval. More explicitly, for a frequency $\psi$, the  sum  over $y_d\leqslant p<y_t$ has
main term 
$$z\log\bigg\{{\lambda_t(1+|\psi|\lambda_d)\over
                    \lambda_d(1+|\psi|\lambda_t)}\bigg\}.$$
This is $z\log(\lambda_t/\lambda_d)+O(1)$ when
$|\psi|\ll_t1/\lambda_t$, and $O(1)$ when
$|\psi|\geqslant1/\lambda_d$. These alternatives cover every
non-empty interval; if both cutoffs coincide the interval is empty.
The extension of \eqref{estimatesump} from frequencies at most one
to frequencies bounded in terms of $t$ only changes its constant.
Also, grouping the non-zero forms by their first containing span
yields
$$\sum_{p<x}{\varrho(p)U_p(\gth)\over p}
\leqslant z\log_2x+
 z\sum_{1\leqslant j\leqslant d}(c_j-c_{j-1})\log\lambda_j+O(1).$$
Consequently
$$\e^{\sum_{p<x}b_p}
\ll (\log x)^z\prod_{1\leqslant j\leqslant d}\lambda_j^{z(c_j-c_{j-1})}
       \lambda_d^{-z2^t}(\lambda_d/\lambda_t)^{za_t}.$$
The factor $(\log x)^z$ in this formula is the contribution of the zero
form. The coefficient $a_t$ is the mass lost on the intermediate
 interval, not the mass retained there. 
\par
On $0\leqslant s_1\leqslant\cdots\leqslant s_d\leqslant d$, define
$$K_z(\s)=
 \sqrt{\lambda_t\lambda_d}
 \prod_{1\leqslant j\leqslant d}\lambda_j^{z(c_j-c_{j-1})}
 \lambda_d^{-z2^t}(\lambda_d/\lambda_t)^{za_t}.$$
After normalisation by $(\log x)^z$, it remains to bound, for each basis,
$$T\int_{0\leqslant s_1\leqslant\cdots\leqslant s_d\leqslant d}
 {K_z(\s)\e^{-F(\e^{\lambda_d})/2}\over\sqrt{\log(2\lambda_d)}}
 \d\s.$$
\par 
Put $z_0=z_t$. We record a useful radial bound: for every non-negative
measurable function $G$ on~$[1,\log x]$,
$$\int_{0\leqslant s_1\leqslant\cdots\leqslant s_d\leqslant d}
 K_{z_0}(\s)G(\lambda_d)\,\dd\s
\ll G(1)+G(\log x)+\int_1^{\log x} G(\lambda){\dd\lambda\over\lambda}.$$
To verify this, write $u=s_t$, $v=s_d$ and integrate first in
$s_1,\ldots,s_{t-1}$. Since $z_0c_j<j$ for $j<t$, successive
integration furnishes the bound $R(u)^{t-1-z_0c_{t-1}}$ for these
variables. The remaining kernel is therefore at most
$$R(u)^{D-1}/R(v)^{D+1},\qquad
D:=t-z_0c_t+z_0a_t-\dm.
$$
This also applies for $t=1$, when there are no previous variables.
Moreover,
$$D\geqslant\alpha(t+z_0)-\dm>0.$$
For $t=1$ the lower bound is $3/2-\sqrt2$; for $t\geqslant2$ it is
again positive. Thus the integral in $u\in[0,v]$ is $\ll 1/R(v)$.
The intervals $0\leqslant v\leqslant 1/\log x$, $1/\log x<v<1$, and
$1\leqslant v\leqslant d$ give, respectively, $G(\log x)$, the integral term, and $G(1)$ in the
asserted bound, using
$\lambda=1/v$ on the middle interval.
\par
Suppose first that $z=z_0$, so that $F=f_T$. Apply the radial bound
with
$$G(\lambda)={\e^{-f_T(\e^\lambda)/2}\over\sqrt{\log(2\lambda)}}.$$
The endpoint terms are bounded. Put $a=\log T$ and $b=\gb_{t,z_0}>0$.
With $r=\log(\lambda+\log3)$, the integral term is at most a
constant times
$$\int_{\log(1+\log3)}^\infty
 {\e^{-(\gc/2)\min(a+r,(r-ba)^2/a)}\over\sqrt{1+r}}
 \,\dd r\ll1.$$
Explicitly, the exponential is at most
$\e^{-c(a+r)/2}+\e^{-c(r-ba)^2/(2a)}$.
The first term is integrable uniformly in $a\geqslant\log3$.
For the second, split at $r=ba/2$ and $r=2ba$: on the central
interval the denominator is of order $\sqrt a$ and the Gaussian
has width of order $\sqrt a$; the two tails are uniformly bounded
by Gaussian decay. Bounded values of $a$ cause no difficulty.
This proves \eqref{EEM2/taut+1} at $z=z_t$.
\par 
Finally suppose that $z>z_0$ and $F=0$. Write $C=2^t-1$ and
$$\kappa_t=\max\{0,\,2^{t-1}-1,\,(2^t-1)/\sqrt2-\alpha\},
\qquad \eta_t=C-\kappa_t>0.$$
The coefficient of $z$ in $\log K_z$ is
$$H:=\sum_{j<t}c_j\log(\lambda_j/\lambda_{j+1})
 +(c_t-a_t)\log(\lambda_t/\lambda_d)+C\log\lambda_d.$$
Since $c_j\leqslant\kappa_t$ for $j<t$ and
$c_t-a_t\leqslant\kappa_t$, we have
$$H\leqslant C\log_2x-\eta_t\log((\log x)/\lambda_d).$$
Therefore, with $\varepsilon:=\eta_t(z-z_0)>0$,
$$K_z\leqslant (\log x)^{C(z-z_0)}K_{z_0}(\lambda_d/\log x)^\varepsilon.$$
The radial bound, now applied with
$G(\lambda)=(\lambda/\log x)^\varepsilon/\sqrt{\log(2\lambda)}$, has
bounded endpoint terms and integral at most
$$ {1\over\sqrt{\log2}}\int_1^{\log x}
       (\lambda/\log x)^\varepsilon{\dd\lambda\over\lambda}
 \leqslant {1\over\varepsilon\sqrt{\log2}}.$$
Since $C(z-z_0)=(2^t-1)z-t=\beta-z$, this proves
\eqref{EEM2/taut+1*}. The inclusion $E_T^{t,z}\subset E_T$ then provides
\eqref{EEM2/taut+1} also for $z>z_t$.
Summing over the finitely many bases completes the proof.
\smallskip
(iii) We now prove \eqref{mu2tgen}. By \eqref{EEM2/taurgeq} and \eqref{EEM2/taur}, the required estimate holds when $t$ is integral. Consider the case $t<k:=\plaf t$. Let $\alpha:=t/k$, $a:=2^{t-k}$, $b:=2^t$ and define $\varrho_a(n):=a^{\omega(n)}\varrho(n)$ and $\varrho_b$ similarly. Since $a^\alpha b^{1-\alpha}=1$, Hšlder's inequality applied to the unnormalised sums provides
$$\EE_{x,\varrho}\big(\mu_2^t\big)\ll {\{(\log x)^{az}\EE_{x,\varrho_a}(\mu_2^k)\}^\alpha(\log x)^{bz(1-\alpha)} \over (\log x)^z}\cdot\eqdef{Holdun}$$
Note that this interpolation is valid for any positive exponent. Under the extra hypothesis $z>\gz_t:=k2^{k-t}/(2^k-1)$, we have $az>z_k$, so \eqref{EEM2/taur} applies. The resulting exponent of $\log x$ is exactly
$$\alpha(2^kaz-k)+(1-\alpha)bz-z=(2^t-1)z-t=\beta-z.$$
\par 
We now remove the restriction $z>\gz_t$. Recall that we may assume $t>1$. For $p\nmid m$, the exact quadratic identity is
$$\mu_2(mp)=\mu_2(m)+{N_{1,2}(m,p)\over \tau(m)},\quad 0\leqslant {N_{1,2}(m,p)\over \tau(m)}\leqslant \mu_2(m), $$
with notation \eqref{defNjq}. By the mean value theorem
$$\mu_2(mp)^t-\mu_2(m)^t\ll_t {\mu_2(m)^{t-1}N_{1,2}(m,p)\over \tau(m)}.\eqdef{incr}$$
Now, expanding $N_{1,2}$ into $\tau(m)^2$ triangular kernels $(1-|\log p-\log d'/d|)^+$ easily yields
$$\sum_{\sqrt{y}<p\leqslant y}{\varrho(p)N_{1,2}(m,p)\over p}\ll {\tau(m)^2\over \log 3y}\quad(y\geqslant 1).$$
\par 
Apply \ref{lemmerhof} to $f=\mu_2^t$ and the above increment. Using \eqref{incr}, we obtain
$$\eqalign{\EE_{x,\varrho}\big(\mu_2^t\big)&\ll1+\int_1^{x^2}\som y_{m\geqslant 1} {\varrho(m)\tau(m)\mu_2(m)^{t-1}\over m}{\dd y\over y(\log 3y)^{z+2}}\cr&\ll 1+\int_1^{x^2}{\EE_{y,{\varrho_2}}(\mu_2^{t-1})\over y(\log 3y)^{2-z}}\d y,\cr}\eqdef{recEmu2t}$$
with the above notation $\varrho_2(n)=2^{\omega(n)}\varrho(n)$.
\par 
Put $k:=\plaf{t-1}<t$. The function $h(u)=u2^u/(2^u-1)$ is strictly increasing for $u>0$ since $2^u>1+u\log 2$. Consequently, $2z_t>\gz_{t-1}$. Thus \eqref{Holdun} is available whenever $z>z_t$. Substituting into \eqref{recEmu2t} yields
$$\EE_{x,\varrho}\big(\mu_2^t\big)\ll1+\int_1^{x^2}{\dd y\over y(\log 3y)^{z-\beta+1}}\ll (\log x)^{\beta-z},$$
since $\beta-z>0$. 
\qed

 \Propl{lemmeck}{Let  $t\in\N^*$. For any basis $\B$ in $B$, we have
$$c_k(\B)  \leqslant 2^k-1\qquad (1\leqslant k\leqslant t).\eqdef{ineglemmeck}$$}
 \nid Let $V\subset \r^t$ be the $k$-dimensional linear space spanned by $w_1,\ldots,w_k$. Choose a set $J$ of $k$ coordinates indices such that the coordinate projection $\pi_J:V\to\r^J$ is injective. A non-zero $k\times k$ minor of a matrix for a basis of $V$ provides such a choice. For every $w\in V$, we have trivially $|w|\geqslant |\pi_J(w)|$. hence
$$1+c_k(\B)=\sum_{w\in V\cap\W_t}{1\over 2^{|w|}}\leqslant \sum_{v\in\{-1,0,1\}^J}{1\over 2^{|v|}}=\Big(1+\dm+\dm\Big)^k=2^k.$$
Injectivity ensures that no term is counted twice and that only $w=0$ projects to 0.\qed
  \medskip

\paraun{Bounding integral moments for bounded $q$}

Recall definitions \eqref{defMq} and \eqref{defmuq}. In  section \ref{prgmoms} below, we generalize the method displayed in  \citer{KT23}, which is based on bounding $\EE_{x,\1}(\1_S\mu_q)$ for sets $S$ of decreasing size and  applying Markov's inequality to estimate  $\PP_{x,\1}(\mu_q>T)$. The approach may easily  be adapted in order to include a weight $\varrho$. However, the initialisation step requires a bound for  $\PP_{x,\varrho}(\mu_q> T)$  decreasing as $1/T^t$. This is provided by estimating $\EE_{x,\varrho}(\1_S\mu_q^t)$.  The generality introduced in \ref{lemmerhof} then enables an iteration of the process.
 \par 
This section is devoted to establishing this initialization step in the case of integral $t$. In the next section, the result is extended to real $t$.
\par 
We now describe the sets $S$ needed in further estimates for conditional expectations.\par 
Let $\gamma\in\{0,1\}$ and  $\{\vartheta_{j,T}^\gamma\}_{j=0}^{\infty}\in(\r^+)^{\N}$ denote a sequence such that $ \vartheta_{0,T}^\gamma= \vartheta_{1,T}^\gamma=1$.  Define further
$$\leqalignno{H_{T,\gamma}^{q}&:=\big\{n\in E_T: M_j(n)\leqslant \tau(n)\vartheta_{j,T}\ (1\leqslant j\leqslant q)\big\},&\eqdef{HTq}\cr
H_{T,\gamma}^{t,z,q}&:=\big\{n\in E_T^{t,z}: M_j(n)\leqslant \tau(n)\vartheta_{j,T}\ (1\leqslant j\leqslant q)\big\},&\eqdef{HTqg}
\cr
H_T^{t,z,q}&:=H_{T,\delta}^{t,z,q},&\eqdef{HTztq}\cr}$$  
so that $H_{T,\gamma}^{0}=H_{T,\gamma}^{1}= E_T $. Note that, by \eqref{recMq}  {\it infra}, $M_j(n)/\tau(n)$ is multiplicatively increasing, so if $n\in H_T^{t,z,q}$ then any divisor of $n$ also lies in~$H_T^{t,z,q}$.
\par 
The quantities  $\vartheta_{j,T}$ will be formally defined later---see \eqref{defthqT}.   For now, we only specify that
$$ \vartheta_{0,T}=\vartheta_{1,T}=1,\qquad \vartheta_{2,T}\asymp T (\log T)^{\gamma},\qquad \vartheta_{j,T}\ll_jT^{j-1 }(\log T)^{\gamma (j-1) }\quad(j\geqslant 3).\eqdef{majthqT}$$
\par
When conditioning over $H_{T,\gamma}^{q}$, we shall need estimates with $\gamma=1$ when $\delta=1.$  
 \par 
Recall notation $\som x$  and \eqref{defmuq}, and put 
 $$\T_{t,z}(x;q):=\PP_{x,\varrho}\big(H_{T}^{t,z,q-1}\big)\EE_{x,\varrho}\Big( \mu_q \Big|H_{T}^{t,z,q-1}\Big)=\prod_{p<x}\Big({1\over 1+\varrho(p)/p}\Big)\som{x}_{n\in H_{T}^{t,z,q-1}}{\varrho(n)\mu_q(n)\over n}\cdot\eqdef{Sqx}$$
\par 
 \ref{Sq} below furnishes a recursive bound for $\T_{t,z}(x;q)$. The initialisation step described above requires a more general estimate that will only be used for bounded values of $q$.  
\par 
  We thus consider the quantities 
$$ 
\leqalignno{\qquad \gS^\gamma_{r,z}(x;q)&:=\PP_{x,\varrho}\big(H_{T,\gamma}^{q-1}\big)\EE_{x,\varrho}\bigg( \mu_q^r \Big|H_{T,\gamma}^{q-1}\bigg)=\prod_{p<x}\Big({1\over 1+\varrho(p)/p}\Big)\som{x}_{n\in H_{T,\gamma}^{q-1}}{\varrho(n)\mu_q(n)^r \over n}\cdot&\eqdef{Sqxrt}}$$ 
Note that the conditioning is over $H_{T,\gamma}^{q-1}$: we shall exploit the fact that this does not depend on $t$ nor $z$.
\Propp{Sbfq}{Let $t\in\N^*,\, r\in\N^*\cap[1,t],\,z\in\r^+,\,\varrho\in\M_z,\,q\geqslant 3$, $\gamma\in \{0,1\}$, and assume  
 $$\cases{Ê z\geqslant z_t &  if $\gamma=1$,\cr z>z_t &  if $\gamma=0$.}\eqdef{condisurgamma}$$    
The sets $H_{T,\gamma}^{j}$ being defined by  \eqref{HTq} with $\theta_{j,T}$ satisfying \eqref{majthqT},  we have, uniformly for $T\geqslant 2$, $x\geqslant 2$, 
$$\gS^\gamma_{r,z}(x;q)\ll T^{r(q-2)}(\log T)^{\gamma r(q-1)}(\log x)^{\beta-z}.\eqdef{majSqxt-int}$$}
\nid 
Put
$$\leqalignno{
N_{j,q}(n,v)&:=\int_\r\Delta(n,u)^j\Delta(n,u-\log v)^{q-j}\d u\quad (n\geqslant 1,\,v\geqslant 1,\,1\leqslant j\leqslant q)&\eqdef{defNjq}\cr 
W_q(n,v)&:=\sum_{1\leqslant j\leqslant   q/2}{q\choose j}N_{j,q}(n,v)\leqslant  2^{q} M_q(n)\quad (n\geqslant 1,\,v\geqslant 1).&\eqdef{Wqnv}
}$$
From \citeplus{Te85}{(17)}, we have
$$2M_q(n)\leqslant M_q(np)\leqslant 2^qM_q(n),\qquad  M_q(np)\leqslant 2M_q(n)+2W_q(n,p)\quad(p\nmid n),\eqdef{recMq}$$
and so, with notation \eqref{defmuq},
$$\mu_q(n)\leqslant \mu_q(np)\leqslant 2^{q-1}\mu_q(n),\quad \mu_q(np)\leqslant  \mu_q(n)+{W_q(n,p)\over \tau(n)}\quad(p\nmid n).\eqdef{rec2Mq}$$
For $q\geqslant 2$ and any $v\geqslant 1$, we also have
$$W_q(n,  v)\leqslant 2^{q-1}\int_\r\Delta(n,u)\Delta(n,u-\log v)\big\{ \Delta(n,u-\log v)^{q-2}+\Delta(n,u)^{q-2}\big\}\d u.\eqdef{inegWq}$$
Define $f:=\1_{H_{T,\gamma}^{q-1}} \mu_q^r  $. By \eqref{rec2Mq} and the inequality $(a+b)^r\leqslant a^r+O(ba^{r-1})$ valid for $b\ll a$, we have, for a suitable constant  $C_{q,r}$, 
 $$f(np) \leqslant f(n)+C_{q,r}{W_q(n,p)\over \tau(n)}\1_{H_{T,\gamma}^{q-1}}(n) \mu_q( n)^{r-1}\cdot $$ 
 \par 
\ref{lemmerhof} is applicable to $\varrho\in\M_z$  and $w$    defined  by $$
w(m,p)= {W_{q,r}(m,p):=C_{q,r}}{W_q(m,p)\over \tau(m)}\1_{H_{T,\gamma}^{q-1}}(m) \mu_q(m)^{r-1} .$$
 Moreover, appealing to \eqref{defNjq} and \eqref{inegWq},   we get,  for all $m\geqslant 1$, 
 $$\eqalign{\int_\r W_{q,r}(m, \e^v)\d v
&\ll   ÊÊ   {M_{q-1}(m)M_{q}(m)^{r-1} \over \tau( m )^{r-1}}\cdot \cr}\eqdef{intWbfq}
$$ 
\par
It follows that
$$\eqalign{\gS^\gamma_{r,z}(x;q)& \ll 1+  
\int_1^{x^2} \som y_{m\in  H_{T,\gamma}^{q-1}}{\varrho(m)\over m  }\sum_{\sqrt{y}\leqslant p< y}{W_{q,r}(m,p)\over p}  {\d y \over y (\log 3y)^{1+z} },}$$
 whence, by partial integration---see \eqref{somWnp} {\it infra}---, 
 $$\eqalign{\gS^\gamma_{r,z}(x;q)
& \ll 1+ 
\int_1^{x^2}{ Z^\gamma_{r,z}(y;q) \over y (\log 3y)^{2+z} }\d y,}\eqdef{recgSq}$$
where we have set
$$Z_{r,z}^{\gamma}(y;q):= \som y_{m\in  H_{T,\gamma}^{q-1}}{\varrho(m)\over m }\int_\r  W_{q,r}(m, \e^v)\d v
\ll
\som y_{m\in  H_{T,\gamma}^{q-1}}{\varrho(m)M_{q-1}(m)\mu_q(m)^{r-1}\over m} \cdot\eqdef{defZbfq}$$
\par 
We now set out to proving \eqref{majSqxt-int} by induction on  $r\geqslant 1$.\par 
\par 
We start with the diagonal case $r=t=1$. If $z>1$, start with $q=3$. Since $\tau(m)\leqslant T\log 3y$ on the
support, \eqref{EEM2/taur} gives
$$
 Z^\gamma_{1,z}(y;3)
 \ll T(\log 3y)\som y_{m\geqslant 1}\frac{\mu(m)^2\varrho(m)\mu_2(m)}{m}
 \ll T(\log 3y)^{2z}.
$$
This implies $\gS^\gamma_{1,z}(x;3)\ll T (\log 3x)^{z-1}$. For $q>3$, nesting of the sets furnishes
$$
 Z^\gamma_{1,z}(y;q)
 \ll T(\log 3y)^{z+1}\gS^\gamma_{1,z}(y;q-1),
$$
and induction yields $\gS^\gamma_{1,z}(x;q)\ll T^{q-2}(\log 3x)^{z-1}$.
This is sufficient for either value of $\gamma$.
\par 
 For $z=1$, only $\gamma=1$ is required. On the support,  $$M_{q-1}\leqslant \vartheta_{q-1,T}\tau,\quad\tau\leqslant \min\Big(\tau,\sqrt{T(\log 3y)\tau}\Big).$$ Hence
$$Z^1_{1,1}(y;q)
 \ll \vartheta_{q-1,T}\min\Big\{(\log 3y)^2,\sqrt T (\log 3y)^{\sqrt2+1/2}\Big\}.
$$
As $s:=3/2-\sqrt2>0$, integration in \eqref{recgSq} yields
$$
 \gS^1_{1,1}(x;q)
 \ll1+\vartheta_{q-1,T}\int_{\log3}^\infty
       \min\{1,\sqrt T/v^{s}\}{\dd v\over v}
 \ll \vartheta_{q-1,T}\log T
 \ll T^{q-2}(\log T)^{q-1}.
$$
\par 
We next perform induction on $t\geqslant 2$, still assuming $r=t$.\par 
Suppose first $z>z_t$. Write $\varrho_2=2^{\omega}\varrho$. The cutoff on $M_{q-1}$ in \eqref{defZbfq} gives
 $$Z^\gamma_{t,z}(y;q)
 \ll T^{q-2}(\log T)^{\gamma(q-2)}(\log 3y)^{2z}
             \gS^\gamma_{t-1,2z}(y;q).
$$
Since $2z_t>z_{t-1}$, the induction hypothesis applies and yields
$$
 Z^\gamma_{t,z}(y;q)
 \ll T^{t(q-2)}(\log T)^{\gamma\{t(q-1)-1\}}(\log 3y)^{\beta(t,z)+1}. \eqdef{eq:ZA}
$$
If $z>z_t$, \eqref{majSqxt-int} for $r=t$ follows by integration.
\par 
If $z=z_t$ and $\gamma=1$, \eqref{eq:ZA} still holds, because the
tilted parameter $2z_t$ is strictly supercritical for $t-1$. Set
$$
 b_t=\log(2z_t/z_{t-1})/\log2>0,\quad
 \varrho_1=\varrho(z_{t-1}/z_t)^\omega,\quad
 d_t=1+z_t-z_{t-1}-b_t>0.
$$
To check the last inequality, use the increasing function
$F(u)=u/(1-1/2^{u})$ with
$$a=F(t-1)<b=F(t),\quad
 -d_t=a-b-\frac{\log(a/b)}{\log2}<0,
$$
because $a,b>1/\log2$ and $a\mapsto a-\log a/\log2$ is increasing
on that interval. Since $\varrho(m)\tau(m)=\varrho_1(m)\tau(m)^{b_t}$,
the critical induction hypothesis at $(t-1,z_{t-1})$ furnishes
$$
 Z^1_{t,z_t}(y;q)
 \ll T^{t(q-2)+b_t}(\log T)^{t(q-1)-1}(\log 3y)^{z_{t-1}+b_t}.
 \eqdef{eq:ZB}
$$
Split the integral in \eqref{recgSq} at $(\log 3y)=T^{b_t/d_t}$,
using \eqref{eq:ZA} below and \eqref{eq:ZB} above this point.
The first integral costs $\ll\log T$; the second costs $\ll1$
after the factor $T^{b_t}$ is canceled. This completes the proof of the diagonal case.
\par 
Hšlder's inequality for the subprobability measure restricted to
$H_{T,\gamma}^{q-1}$ furnishes
$$
 \gS^\gamma_{r,z}(x;q)
 \leqslant \{\gS^\gamma_{t,z}(x;q)\}^{r/t}
 \ll T^{r(q-2)}(\log T)^{\gamma r(q-1)}(\log 3x)^{(\beta-z)r/t}.
$$
Since $\beta-z\geq0$, this implies  \eqref{majSqxt-int}. Bounded $x$ is absorbed
into the constant. \par 
This completes the proof of  \ref{Sbfq}. \qed

\medskip
 \bigskip
 \paraun{Bounding real moments for bounded $q$}
Recall definition \eqref{Sqxrt}. 
When the parameter $t$ is integral, \ref{Sbfq} enables  initialising the induction conducted in the proof of \ref{Sq}  below  by providing  the upper bound \eqref{majSqxt-int} for $\gS_{t,z}(x;q)$.
 The next statement extends this to real values of $t$ provided $$z\geqslant \gz_t:=2^{\plaf t -t}\plaf t /(2^{\plaf t }-1).$$ 
 Note that $\gz_t=z_t$ if $t\in\N^*$.
  \Propl{EMqt/taut-gen}{ Let $z>0$, $t>0$, $\varrho\in\M_z$.\par 
 {\rm (i)} Let $q\geqslant 3$, $\gamma\in \{0,1\}$. Assume  $z > \gz_t$ if $\gamma=0$, and $z\geqslant \gz_t$ if $\gamma=1$.  Then we have $$\gS^\gamma_{t,z}(x;q)\ll T^{t(q-2)} (\log T)^{\gamma t(q-1) }(\log x)^{\beta-z}\qquad (x\geqslant 2,\,T\geqslant 3).\eqdef{faEMqt/taut-gen}$$
  \par 
   {\rm(ii)}  Under condition  $z> \gz_t$, we have furthermore
 $$\gS^\gamma_{t,z}(x;2)\ll  (\log x)^{\beta-z}\qquad (x\geqslant 2,\,T\geqslant 3).\eqdef{faEMqt/taut-genq=2}$$ \par 
{\rm(iii)}  If $z\geqslant \gz_t$, the bound \eqref{EEM2/taut+1} holds. If $z>\gz_t$ the bound \eqref{EEM2/taut+1*} is satisfied, and, more accurately, with $c_t:=\plaf t-t$,
$$\EE_{x,\varrho}\big(\1_{E_T^{t,z}}  \mu_2^{t+1}\big)  \ll T (\log x)^{\beta-z}\e^{-c_tf_T(x)/(\plaf t+1)}
\qquad ( x\geqslant 3 ).\eqdef{lastone}$$
  }

 \rem
 When $t$ is not integral, the assumption $z\geqslant \gz_t$ implies $z>z_t$.  
\goodbreak\smallskip
  \nid  (i) As noticed above, when $t\in\N^*$, estimate \eqref{faEMqt/taut-gen} follows from \eqref{majSqxt-int}. 
  \par 
 Assume then $t\in\r^*\sset\N$, and put $k:=\plaf t>t$. 
Let $u,\,v$ be positive real numbers such that $uv=z$. By Hšlder's inequality and~\eqref{majSqxt-int} with $r=k$, we have, provided $u^{k/t}\geqslant \gz_k=k/(2^k-1)$,
  $$\eqalign{\gS^\gamma_{t,z}(x;q)& \leqslant \EE_{x,\varrho}\bigg( { u^{k\omega/t}\mu_q^k\over z^{\omega}}\1_{H_{T,\gamma}^{q-1}}\bigg)^{t/k}\EE_{x,\varrho}\Big({v^{k\omega/(k-t)}\over z^\omega}\Big)^{1-t/k}\ll T^{t(q-2)}(\log x)^{B-z}}(\log T)^{\gamma t(q-1)}
  $$
   with $$B:=\Big(2^ku^{k/t}-k\Big) \frac tk+  v^{k/(k-t)}\Big(1-\frac tk\Big).$$     \par
 Select $u:=2^{-t(1-t/k)} z^{t/k}$ so that $B=z2^t -t=\beta$. Observe that $u^{k/t}=2^{t-k}z$. Hence condition  $u^{k/t}\geqslant \gz_k$  amounts to
  $z\geqslant  2^{k -t}k/(2^{k }-1)=\gz_t,$ and $u^{k/t}>\gz_k$  amounts to  $z>\gz_t$. This furnishes \eqref{faEMqt/taut-gen} under the stated hypotheses. 
     \par\goodbreak
  (ii) Similarly, \eqref{faEMqt/taut-genq=2} follows from \eqref{EEM2/taur} on selecting the same value for $u$. 
   \par
  (iii) Consider first the validity of \eqref{EEM2/taut+1} under the more general hypothesis $z\geqslant \gz_t$. Introducing as previously parameters $u,\,v$ such that $z=uv$. Writing $u_1:=u^{(k+1)/(t+1)}$, $v_1:=v^{(k+1)/(k-t)}$ for legibility, we have, appealing to the bound $2^{\omega(n)}\leqslant T\log 3x$ $(n\in E_{T},\, P^+(n)<x)$ and using \eqref{EEM2/taut+1} for exponent $k+1$,
  $$\eqalign{\EE_{x,\varrho}\bigg(\1_{E_T^{t,z}} \mu_2^{t+1}\bigg) & \leqslant \EE_{x,\varrho}\bigg(\1_{E_T^{t,z}} {\mu_2^{k+1} u_1^{\omega}\over z^{\omega}}\bigg)^{(t+1)/(k+1)}\EE_{x,\varrho}\Big(\1_{E_T^{t,z}}{v_1^{\omega}\over z^\omega}\Big)^{(k-t)/(k+1)}\cr
  & \leqslant \EE_{x,\varrho}\bigg( \1_{E_T^{t,z}}{\mu_2^{k+1} u_1^{\omega}\over z^{\omega}}\bigg)^{(t+1)/(k+1)}(T\log x)^{(k-t)/(k+1)}\EE_{x,\varrho}\Big({v_1^{\omega}\over (2z)^\omega}\Big)^{(k-t)/(k+1)} 
  } 
  $$
Selecting $u:=2^{- (k-t)(t+1)/(k+1)} z^{(t+1)/(k+1)}$  and checking that $2^{k}u_1=2^{t}z$ so that  $E_T^{t,z}=E_T^{k,u_1},$ we may conclude as previously.  The proof that \eqref{EEM2/taut+1*} is valid under the extended hypothesis $z>\gz_t$ is similar. Estimate \eqref{lastone} is derived by exploiting the condition $2^{\omega(n)}\leqslant \e^{-f_T(x)}T \log 3x$, valid for $n\in E_{T}^{t,z},\, P^+(n)<x$.
   \qed 
 \smallskip
  \Propl{lemmereg}{Let $z>0$, $t\geqslant 1$, $\varrho\in\M_z$.
  \par 
 {\rm(i)}  Let $q\geqslant 3$, $\gamma\in \{0,1\}$. Assume  $z > z_t$ if $\gamma=0$, and $z\geqslant z_t$ if  $\gamma=1$. Assume furthermore that \eqref{condisurgamma} holds. Then \eqref{faEMqt/taut-gen} persists.
 \par {\rm(ii)}  
 If $z>z_t$ the bounds \eqref{EEM2/taut+1}  and \eqref{EEM2/taut+1*} are satisfied. Moreover, 
 if $z\geqslant z_t$, we have
  $$  \EE_{x,\varrho}\big(\1_{E_T^{t,z}}  \mu_2^{t+1}\big)\ll T(\log T)^{\delta/2} (\log x)^{\beta-z}.\eqdef{4.5plus}$$
}\smallskip
 \nid  (i) We may plainly assume that $t$ is not an integer and hence that $t>1$. 
Arguing as for \eqref{recgSq}, we get
 $$\gS^\gamma_{t,z}(x;q)\ll 1+\int_1^{x^2}{ Z_{t,z}^\gamma(y;q)\over y  (\log 3y)^{2+z}  }\d y.\eqdef{ZL25}$$  
 Appealing to the inequality $M_{q-1}(m)\leqslant \theta_{q-1,T}\tau(m)$ taking \eqref{majthqT} into account, we infer that
 $$Z^{\gamma}_{t,z}(y;q)\ll T^{q-2}(\log T)^{\gamma(q-2)}(\log 3y)^{2z}\gS^\gamma_{t-1,2z}(y;q).\eqdef{ereg}$$ 
Now, observe that $2z_t > \gz_{t-1}$ since $t >  \plaf{t-1}$.\note{This follows from the fact that $s\mapsto s2^s/(2^s-1)$ is increasing.}
  We may hence deduce from \eqref{faEMqt/taut-gen} that
 $$\gS^\gamma_{t-1,2z}(y;q)\ll T^{(q-2)(t-1)}(\log T)^{(t-1)(q-1)\gamma}(\log 3y)^{\beta+1-2z}, \eqdef{maj1Zj*} $$ 
 since $ \beta(t-1,2z)=\beta(t,z)+1=\beta+1$. 
 Indeed, since $2z>z_{t-1}$, the bound for $\gS_{t-1,2z}$   is available whenever  $\gamma\in \{0,1\}$. 
 Carrying back into~\eqref{ereg}, we obtain 
  $$Z^{\gamma}_{t,z}(y;q)\ll T^{t(q-2)}(\log T)^{ \gamma\{t(q-1)-1\}} (\log 3y)^{\beta+1},$$ which yields the required estimate \eqref{faEMqt/taut-gen} by \eqref{ZL25} provided $\delta=0$, and so $\beta+1-(z+2)>-1$.\par 
  If $\delta=1$, \ie $z=z_t$  (hence $\gamma=1$), this is insufficient to get the expected bound.
  We circumvent the difficulty by introducing a parameter  $\gz_{t-1}^*\in\big]\gz_{t-1},2z_{t}\big[$  and writing
 $$\eqalign{Z^{1}_{t,z}(y;q)&\ll\vartheta_{q-1,T}\som y_{m\in H_{T,  1}^{q-1}}{\varrho_2(m) 
 2^{\lambda_t\omega(m)} \mu_{q}(m)^{t-1}\over m}\cr
 &\ll \vartheta_{q-1,T}\big\{ T(\log 3y)\big\}^{\lambda_t}\som y_{m\in H_{T,  1}^{q-1}}{\varrho_2(m) 
 \mu_{q}(m)^{t-1}\over m}
\cr}$$
 with $\varrho_2(m):=\varrho(m)(\gz_{t-1}^*/z_{t})^{\omega(m)}$ and $\lambda_t:=\{\log (2z_{t}/\gz_{t-1}^*)\}/\log 2$. Taking \eqref{majthqT} into account,  the bound \eqref{faEMqt/taut-gen} yields 
 $$Z^1_{t,z} (y;q)\ll T^{t(q-2)+\lambda_t}(\log T)^{t(q-1)\gamma-\gamma}(\log 3y)^{\mu_t}, \eqdef{maj2Zj*}$$ with 
 $$\mu_t:=\beta(t-1,\gz_{t-1}^*)+\lambda_t <z_t+1,$$
 provided $\gz_{t-1}^*\in]\gz_{t-1},2z_t[$ is chosen sufficiently close to $\gz_{t-1}$. 
 To see the latter inequality, observe that, taking limits, $\mu_t-z_t-1=\varphi(a,b):=(a-b)-\{\log (a/b)\}/\log 2$ with $a=f(\plaf{t-1})$,  $b=f(t)$, and $f(s):=s/(1-1/2^s)> 1/\log 2$. Since $\varphi(a,b)$ is a non-decreasing function of $a$ and   $1/\log 2<a< b$, the inequality follows.\par
  \goodbreak
Now split the integral in \eqref{ZL25} at $\log 3y=T^\nu$ with $\nu:= \lambda_t /(1+z_t-\mu_t)$, 
 applying \eqref{maj2Zj*} to the upper range and  \eqref{maj1Zj*}  to the lower range.
 This yields \eqref{faEMqt/taut-gen} as required.
 \par 
 \smallskip 
(ii)   Let us focus on proving \eqref{4.5plus} for real $t\geqslant 1$, $z\geqslant z_t$, and leave  the proof of \eqref{EEM2/taut+1*} for $z>z_t$ to the reader. We may plainly assume $t\notin\N^*$.
 Arguing as in the proof of \eqref{recgSq}, we get
 $$\EE_{x,\varrho}\Big(\1_{E_T^{t,z}}\mu_2^{t+1}\Big)\ll 1+\int_1^{x^2}{ Z_{t,z} (y;T)\over y  (\log 3y)^{2+z}  }\d y, $$
 with now 
 $$Z_{t,z} (y;T):=\som y_{m\in E_{T}^{t,z}}{\varrho(m)\tau(m) \mu_2(m)^{t}\over m }\cdot$$
Since $2z_t > \gz_{t-1}$, we may apply \eqref{lastone}  with $(2z, t-1)$ in place of $(z,t)$ to get
 $$Z_{t,z} (y;T)\ll T(\log 3y)^{1+\beta}\e^{-(\plaf t -t)f_T(y)/(\plaf t+1)}. $$
This implies 
 $$  \EE_{x,\varrho}\big(\1_{E_T^{t,z}} \mu_2^{t+1}\big)\ll T(\log T)^{\delta/2} (\log x)^{\beta-z}.$$
 \vskip-5mm \qed 
\bigskip 
\medskip
  \paraun{Bounding moments inductively}\drefun{moms}

Throughout this section, $t\geqslant1$ and $z\geqslant z_t$ are fixed,
$\beta=2^tz-t$, $\delta=\1_{z=z_t}$, and $T\geqslant3$.  Recall definition \eqref{Sqx} for $\T_{t,z}(x;q)$.\par 
 In the sequel of this section, all constants  may depend on $t,z$, the constants defining $\M_z$, and the
fixed constant in $f_T$, but not on $x,T,q$. Constants in the finite-rank
initialisation may additionally depend on the rank.
\par 
Given a constant $A\geqslant 1$, we define 
$$\vartheta_{j,T}:={j!\over j^2}A^{j -1}T^{j-1}(\log T)^{\delta(j-1)} \quad(j\geqslant 1),\eqdef{defthqT}$$
so that $\vartheta_{1,T}=1$. It can be easily checked that we have
$$\sum_{ 1\leqslant j\leqslant   q/2}{q\choose j}\vartheta_{j,T}\vartheta_{q-j,T}\leqslant {2\pi^2\vartheta_{q,T} \over 3A T(\log T)^{\delta/2}  }\quad(q\geqslant 3).\eqdef{inegthqT}
$$  
\Propp{Sq}{Let  $t\geqslant 1,\,z\geqslant z_t$, $\varrho\in\M_z$, and let $A$ be a sufficiently large  constant depending on the fixed parameters. With the above notation, we have 
$$\T_{t,z}(x;q)\leqslant
 {\vartheta_{q,T}\over q^2T^t}(\log3x)^{\beta-z}
 \qquad(x\geqslant2,\ q\geqslant q_t:=\plaf t+1).
\eqdef{majSq}$$}

\nid
 We give the finite-rank initialisation separately, since bounds for
probabilities of exit do not by themselves initialise this first-moment
estimate. 
 \smallskip\noindent{\it Prime averaging and fixed-rank recursion.}
By \eqref{recMq}, the function $\1_{ H_{T}^{t,z,q-1}} \mu_q $
satisfies \eqref{relf} with  $w(n,p)=\1_{ H_{T}^{t,z,q-1}}(n)W_q(n,p)/\tau(n)$.
By partial integration and appeal to the Brun-Titchmarsh inequality as in \citeplus{KT23}{(6.9)},  we get
$$\sum_{\sqrt{y}\leqslant p< y}{\varrho(p)W_q(m,p)\over p}\ll \sum_{ 1\leqslant j\leqslant   q/2}\Big(1+2^j{\log 3y\over y^{1/4}}\Big){q\choose j}{M_{q-j}(m)M_j(m)\over\log y}
 \cdot\eqdef{somWnp}
$$ Applying \ref{lemmerhof}, we therefore obtain
 $$ \T_{t,z}(x;q)\ll 1+\sum_{ 1\leqslant j\leqslant   q/2}{q\choose j}
\int_1^{x^2}\som y_{m\in H_T^{t,z,q-1}}\Big(1+2^j{\log 3y\over y^{1/4}}\Big){\varrho(m)M_{q-j}(m)M_j(m)\over m\tau(m)y {(\log 3y)^{2+z}}}{\d y  }.\eqdef{rec}$$

Let $S$ be a divisor-closed subset of $E_T$, and let
$\sigma\in\M_v$. Put
$$F_{k,v}(y;S):=\EE_{y,\sigma}(\1_S\mu_k),\qquad
 \sigma_2(n):=2^{\omega(n)}\sigma(n).$$
The weight underlying $F$ is understood; replacing $v$ by $2v$ in a
recursion means replacing $\sigma$ by $\sigma_2$.
Log-convexity of the moments implies $M_jM_{k-j}\leqslant\tau M_{k-1}$.
Thus \ref{lemmerhof} and \eqref{somWnp}, for each fixed $k\geqslant2$,
yield
$$F_{k,v}(x;S)\ll_k1+
 \int_1^{x^2}(\log 3y)^{v-2}F_{k-1,2v}(y;S){\dd y\over y}\cdot
\eqdef{initrec}$$
Here $\mu_1=1$. Alternatively, if $S\subset E_T^{t,z}$, retaining the
original weight and using $\tau(m)\leqslant T(\log 3y) \e^{-f_T(y)}$ gives
$$F_{k,z}(x;S)\ll_k1+
 T\int_1^{x^2}\e^{-f_T(y)}F_{k-1,z}(y;S){\dd y\over y\log 3y}\cdot
\eqdef{initrecT}$$
Both inequalities are used only at bounded ranks during initialisation.

 \smallskip\noindent{\it Initialisation at $q_t$.}
Suppose first $t>1$, and set
$$k=\plaf t,\quad r=t-k+1\in]0,1],\quad a=k-t=1-r,
 \quad B=1+\vartheta_{2,T}^a.$$
On any of our sets $S\subset H_T^{t,z,2}$, we have
$\mu_2\leqslant\vartheta_{2,T}^a\mu_2^r$.
For $0<r\leqslant1$ and $v>\gz_r=2^{1-r}$, concavity and a change
of weight give
$$\EE_{y,\sigma}(\mu_2^r)
 \ll (\log 3y)^{(2^r-1)v-r}.
\eqdef{initaux}$$
For completeness, set $c=2^{r-1}$ and $d=2^r$. Since $c^rd^{1-r}=1$,
H\"older's inequality on the unnormalised sums, followed by the integral
quadratic estimate at parameter $cv>1$, gives the exponent
$$r(2cv-1)+(1-r)dv-v=(2^r-1)v-r.$$
This proves \eqref{initaux}; for $r=1$ it is just the quadratic estimate.

At $z=z_t$, the condition $2^{k-1}z>\gz_r$ is equivalent to
$t/(1-2^{-t})>2$, and is therefore strict when $t>1$.
Starting with \eqref{initaux} at weight $2^{(k-1)\omega}\varrho$, apply
\eqref{initrec} successively up to degree $k$ and weight $2^\omega\varrho$.
All intermediate powers are strictly positive. Indeed, with
$h(u)=u/(1-2^{-u})$, the inequalities
$2^{t-s}z_t>z_s$ for $0<s<t$ follow from $h(t)>h(s)$.
We obtain, uniformly for $S\subset H_T^{t,z,2}$,
$$F_{k,2z}(y;S)\ll
 B (\log 3y)^{(2^t-2)z-t+1}.
\eqdef{initpenult}$$
The same estimate holds after replacing $z$ by any sufficiently close
smaller parameter $\xi$, with the set $S$ and its thresholds held fixed.
In particular, if
$$Z_z(y;S):=\som y_{m\in S}{\varrho(m)\tau(m)\mu_k(m)\over m},$$
then $Z_z(y;S)\ll B (\log 3y)^{2^tz-t+1}$.
The last use of the functional lemma therefore proves, for $z>z_t$,
$$F_{q_t,z}(x;S)\ll B (\log x)^{\beta-z}.
\eqdef{initstrict}$$

If $z=z_t$, choose $0<\xi<z$ sufficiently close to $z$, and put
$$b:={\log(z/\xi)\over\log2},\qquad
 d:=z+t-2^t\xi-b>0.$$
Such a choice is possible because $2^tz\log2>1$; it also preserves all
strict conditions preceding \eqref{initpenult}.
For $\varrho_\xi=(\xi/z)^\omega\varrho$, we have
$\varrho=\varrho_\xi\tau^b$. The bound $\tau(m)\leqslant T\log 3y$ on $S$
and \eqref{initpenult} at parameter $\xi$ imply
$$Z_z(y;S)\ll B T^b (\log 3y)^{z+1-d}.$$
Combining this with $Z_z(y;S)\ll B (\log 3y)^{z+1}$, we may bound the final integral by
$$B\int_{\log3}^{\infty}\min\{1,T^bv^{-d}\}{\dd v\over v}
 \ll B\log T.$$
Thus $F_{q_t,z_t}(x;S)\ll B\log T$.
\par 
When $t=1$, take $q_t=2$, $a=0$ and $B=1$. For $z>1$ use the
unrestricted quadratic estimate. For $z=1$ apply the functional lemma
directly to $\mu_2$: its numerator is bounded by $(\log 3y)^2$ and, after the
same change of weight, by $T^b(\log 3y)^{2\xi+b}$, where
$b=\log(1/\xi)/\log2$ and $2-2\xi-b>0$ for $\xi<1$ sufficiently close
to $1$. The preceding integral then gives $F_{2,1}(x;E_T)\ll\log T$.

 \smallskip\noindent{\it The remaining bounded ranks.}
Apply \eqref{initrecT} successively, using the nesting of the sets $H$.
In the strict case the integral of $(\log 3y)^{\beta-z-1}\dd y/y$ is
$\ll (\log x)^{\beta-z}$. At the critical point,
$$\int_1^\infty \e^{-f_T(y)}{\dd y\over y\log 3y}\ll\sqrt{\log T}.
\eqdef{radialT}$$
To check this, put $u=\log_2 3y$ and $v=\log T$. From the definition of
$f_T$, its exponential is at most the sum of
$ \e^{-\gc(v+u)}$ and
$ \e^{-\gc(u-\gb v)^2/v}$, whose integrals are respectively
$\ll1$ and $\ll\sqrt v$.
In particular, for each bounded integer $q\geqslant q_t$,
$$\T_{t,z}(x;q)\ll_q (1+A^a)T^{q-1-t}(\log x)^{\beta-z}
 (\log T)^{\delta(q-t)}.
\eqdef{initfinite}$$
Here $a=\plaf t-t$; the weaker factor $\log T$ per step suffices
in deriving the displayed critical bound. Since $q-t\leqslant q-1$
and $q-1>a$, increasing $A$ establishes \eqref{majSq} simultaneously
for all integers $q_t\leqslant q<2q_t$. The constants in
\eqref{initfinite}, apart from the displayed $1+A^a$, are independent
of $A$: only the bound on $\mu_2$ was used in this initialisation. 
\par 
\smallskip\noindent{\it Induction at the large ranks.}
Let $q\geqslant2q_t$ and assume \eqref{majSq} at all smaller
eligible ranks. In particular, $q-j\geqslant q_t$
whenever $1\leqslant j\leqslant q/2$.

\par 
The contribution of the term involving $y^{-1/4}$ in \eqref{rec} may be handled as in \citer{KT23}, 
using 
the trivial bound 
$$M_{q-j}(m)M_j(m)\leqslant \tau(m)^q\leqslant  T^{q-t-1}(\log 3y)^{q-t-1} \tau(m)^{t+1},$$ valid whenever $m\in E_T$ and $P^+(m)< y$.
Since 
$$\som y_{m\geqslant 1}{\varrho(m)\tau(m)^t\over m}\ll(\log 3y)^{2^tz},$$
this provides an overall term at most
$$ 3^{q}T^{q-t-1}\int_1^\infty{(\log 3y)^{q+\beta-z-2}\over y^{5/4}}\d y\ll C^q q!T^{q-t-1},$$ where~$C$ is a suitable  constant.
 \par 
  Therefore, since $\beta-z\geqslant 0$,
$$ T_{t,z}(x;q)\ll 1+  C^q q!T^{q-t-1} +\sum_{ 1\leqslant j\leqslant   q/2}{q\choose j}
\int_1^{x^2}{Z_{q,j}(y,T)\over y {(\log 3y)^{2+z}} }\d y,$$
with 
$$Z_{q,j}(y,T):=\som y_{m\in H_T^{t,z,q-1}}{\rho(m)\mu_{q-j}(m)M_j(m)\over m}\cdot
$$
\par \goodbreak
We aim at establishing \eqref{majSq} by induction on   $q\geqslant 2q_t$. 
 By \eqref{HTztq}, we have, whenever $m$ is counted in  $Z_{q,j}$,
$$M_j(m)\leqslant \vartheta_{j,T} T (\log 3y)\,\e^{-f_T(y)}  \qquad \Big(m\in H_T^{t,z,q-1},\,P^+(m)<y,\,1\leqslant j\leqslant q-1,\,y\geqslant 1\Big).$$ 
Therefore
$$ Z_{q,j}(y,T)\ll \vartheta_{j,T} T(\log 3y)^{1+z}\,\e^{-f_T(y)}\T_{t,z}(y;q-j).$$ 
For $j\leqslant q/2$, we have $q-j\geqslant q/2\geqslant t$. We may hence bound $\T_{t,z}(y;q-j) $ by the induction hypothesis. This yields
$$ Z_{q,j}(y,T)\ll {\vartheta_{j,T} \vartheta_{q-j,T}\over q^2T^t} T \,\e^{-f_T(y)} (\log 3y)^{1+\beta}
,$$
whence
 $$\eqalign{ \T_{t,z}(x;q)&\ll 1+C^qq!T^{q-t-1}+ \sum_{1\leqslant j\leqslant q/2}{q\choose j}{\vartheta_{q-j,T}\vartheta_{j,T}  \over q^2T^{t-1} }\int_1^{x^2}{ \,\e^{-f_T(y)} (\log 3y)^{\beta-z-1}\over y }\d y  
\cr&\ll 1+C^qq!T^{q-t-1}+ {(\log 3x)^{\beta-z}(\log T)^{\delta/2}\over q^2T^{t-1}}\sum_{1\leqslant j\leqslant q/2}{q\choose j}{\vartheta_{q-j,T}\vartheta_{j,T}}.}\eqdef{majSqavecsomme}$$ 
\par 
By \eqref{defthqT}, $$\eqalign{
\sum_{1\leqslant j\leqslant q/2}{q\choose j}
                     \vartheta_{j,T}\vartheta_{q-j,T}
 &= {\vartheta_{q,T}\over AT(\log T)^\delta}
       \sum_{1\leqslant j\leqslant q/2}{q^2\over j^2(q-j)^2}\leqslant {4\vartheta_{q,T}\over AT(\log T)^\delta}
                    \sum_{j\geqslant1}{1\over j^2}\cdot\cr}$$
                    Since $T\geqslant3$, it follows, with suitable $C$, that
$$
\T_{t,z}(x;q)\leqslant C+C^q q!T^{q-t-1}+{C\vartheta_{q,T}(\log3x)^{\beta-z}\over Aq^2T^t}\cdot
                           $$
We now absorb these three terms separately. First choose $A$ so large
that $C/A\leqslant1/4$. For the second term, the exact quotient is
$${C^q q!T^{q-t-1}\over
      \vartheta_{q,T}(\log3x)^{\beta-z}/(q^2T^t)}
 ={C^q q^4\over
      A^{q-1}(\log T)^{\delta(q-1)}(\log3x)^{\beta-z}}
 \leqslant{C^q q^4\over A^{q-1}}.$$
Increasing $A$ makes this at most $1/4$, simultaneously for all
$q\geqslant2q_t$. Finally,
$${\vartheta_{q,T}(\log3x)^{\beta-z}\over q^2T^t}
 ={A^{q-1}q! T^{q-1-t}
       (\log T)^{\delta(q-1)}(\log3x)^{\beta-z}
\over q^4} \geqslant{A^{q-1}q!\over q^4}\cdot$$
This is at least $4C$ for all $q\geqslant2q_t$ after a further increase
of $A$. The same choice may be made large enough for the previously
established finite-rank initialisation, since its dependence on $A$
is only through $1+A^a$, with $a=\plaf t-t<q-1$.
Thus
$$\T_{t,z}(x;q)\leqslant
 {3\vartheta_{q,T}(\log3x)^{\beta-z}\over 4q^2T^t}
 \leqslant{\vartheta_{q,T}(\log3x)^{\beta-z}\over q^2T^t},$$
which completes the induction, thus implying \eqref{majSq}.
\qed

\bigskip

\paraun{Proof of \ref{thD}}
Recall the sequence $\vartheta_{j,T}$ defined in  \eqref{defthqT}, satisfying \eqref{inegthqT}.
From now on, we consider the sets $H_T^{t,z,q}$ defined in \eqref{HTztq}.
\par
By   
\eqref{EssetET}    and  \eqref{4.5plus}, we derive using Markov's inequality that, for $z\geqslant z_t$  $$\PP_{x,\varrho}\big(E\sset H_T^{t,z,2}\big)\ll {(\log x)^{\beta-z}\over T^t}\qquad (x\geqslant 2,\,T\geqslant 3).\eqdef{PE-HT}$$

\smallskip
Let us put $\lambda:=C_1 T  (\log T)^{\delta} $ and set out to establishing that, for fixed $t\geqslant 1$, $z\geqslant z_t$, we have
$$\PP_{x,\varrho}\Big(\Delta(n)>\lambda\log_2x  \Big)\ll{(\log x)^{\beta-z}
\over T^t}\qquad (x\geqslant 3,\,T\geqslant 1).\eqdef{majPxD}$$
By 
 \eqref{PE-HT}, it is sufficient to show that 
$$\PP_{x,\varrho}\Big(n\in  H_T^{t,z,2},\,\Delta >\lambda\log_2x\Big)\ll{(\log x)^{\beta-z}\over T^t}\cdot$$
However, from \eqref{majSq},   we have, for $q\geqslant q_t$,
$$\PP_{x,\varrho}\Big(H_T^{t,z,q-1}\sset H_T^{t,z,q}\Big)\leqslant \PP_{x,\varrho}\Big(H_{T}^{t,z,q-1}\Big)\EE_{x,\varrho}\Big({\mu_q\over \vartheta_{q,T}}\Big|H_T^{t,z,q-1}\Big)\ll{(\log x)^{\beta-z}
\over q^2T^t}\cdot$$
When $3\leqslant q\leqslant q_t$, we deduce from  \ref{lemmereg} that
$$\eqalign{\PP_{x,\varrho}\Big(H_T^{t,z,q-1}\sset H_T^{t,z,q}\Big)&\leqslant \PP_{x,\varrho}\Big(H_T^{t,z,q-1}\Big)\EE_{x,\varrho}\bigg({\mu_q^t\over \vartheta_{q,T}^t}\Big|H_T^{t,z,q-1}\bigg)\ll{\gS_{t,z}^{\delta}(x;q)\over \vartheta_{q,T}^t}Ê\ll {(\log x)^{\beta-z}
\over q^2T^t}
,\cr}$$
and so  
$$\PP_{x,\varrho}\Big( E\sset \cap_{q}H_T^{t,z,q}\Big)\ll{(\log x)^{\beta-z}
\over T^t}\cdot\eqdef{PE-cupH}$$
\par 
Now observe that, for $n\in\cap_qH_T^{t,z,q}$ and $P^+(n)\leqslant x$, we have
$$\Delta(n)\leqslant 2 M_q(n)^{1/q}\ll \tau(n)^{1/q}\vartheta_{q,T}^{1/q}\ll  q T(\log T)^{\delta }  (\log x)^{1/q}  ,$$
where the first inequality is \citeplus{HT88}{(5.56)}.
Selecting $q=\fl{\log_2x}$ we deduce that, for a suitably large constant $C_1$, we have $\Delta(n)\leqslant \lambda \log_2x$. Thus \eqref{majPxD} follows from \eqref{PE-cupH}.
  \par 
It remains to estimate $\EE_{x,\varrho}(\Delta^t).$  We may assume trivially that $2(\log_2x) \leqslant \Delta(n)\leqslant (\log x)^{z 2^t +t}$.
Indeed, since $\Delta \leqslant \tau$, the contribution of $\Delta(n)>(\log x)^{z2^t+t}$ in $\EE_{x,\varrho}(\Delta^t)$ does not exceed 
$${\EE_{x,\varrho}(\tau^{t+1})\over (\log x)^{z 2^t +t}}\ll {(\log x)^{z2^{t+1}-z}\over (\log x)^{z 2^t +t}} 
\ll  (\log x)^{\beta-z}.$$
 The contribution of those integers $n$ such that $2^j<\Delta(n)/\log_2x \leqslant 2^{j+1}$ is plainly
$$\ll{ 2^{t j} (\log_2x)^t(j+1)^{t\delta}  (\log x)^{\beta-z}\over  2^{t j}}
=  { (\log_2x)^t(j+1)^{t\delta} (\log x)^{\beta-z}
}
 .$$
Summing over $j\ll \log_2x$, 
we finally get
$$\EE_{x,\varrho}(\Delta^t)\ll(\log_2x)^{1+t+t\delta}(\log x)^{\beta-z}.$$
Since, by \ref{2.1}, we have
$$S_{t,\varrho}(x)=\sum_{n\leqslant x}\varrho(n)\Delta(n)^t\ll x(\log x)^{z-1}\EE_{x,\varrho}(\Delta^t),$$ we obtain the upper bound of \eqref{estsDt} as required.
\par 
 To establish \eqref{estsDts}, we apply Hšlder's inequality with parameters $s/t$ and $s/(s- t)$ and invoke~\eqref{estsDt} with $t=s$,  $\delta=1$, yielding
$$\EE_{x,\varrho}\big(\Delta^t\big)\leqslant \Big\{ \EE_{x,\varrho}\big(\Delta^s\big)\Big\}^{t/s}\ll  \big(\log_2x\big)^{2 t+t/s}.$$

\bigskip  
\paraun{Proof of \ref{thDelta3}}
 Write $t=2s$ with $s\geqslant 1$. By \citeplus{HT88}{Exercise 47} (see a proof in \citeplus{BT25}{\S 4}) we have
$$\Delta(n)^2\leqslant 4 M_2(a)M_2(b)\qquad (ab=n),\eqdef{inegDelta2ab}$$
and so, for squarefree $n$,
 $$\Delta(n)^{t}\leqslant  2^t \sum_{ab=n}{M_2(a)^sM_2(b)^s\over 2^{\omega(ab)}}\cdot$$ 
 Appealing to \eqref{mu2tgen} with $\varrho_1(n):=\varrho(n)2^{(s-1)\omega(n)}$,
 we deduce that, provided $2^{s-1}z>z_s$, we have
$$\EE_{x,\varrho}\big(\Delta^{t}\big)\ll {1\over (\log x)^z}\bigg\{(\log x)^{ 2^{s-1}z} \,\EE_{x,\varrho_1}\big(\mu_2^s\big)\bigg\}^2\ll (\log x)^{(2^{t}-1)z-t}.$$
In view of \eqref{Marith-log} and \eqref{triv},  this furnishes the  true order of magnitude.
\bigskip
\medskip
 \paraun{Proof of \ref{thDelta2}}
    Put $\chi_k(n):=\1_{\omega(n)=k}$ $(n\geqslant 1,\,k\geqslant 0)$, so that 
$$\EE_{x,\varrho}(\Delta^2)=\sum_{k\geqslant 0} \EE_{x,\varrho}(\Delta^2\chi_k).\eqdef{decD2}$$

The proof of \eqref{S21x} relies on the following bound.
\Propl{lemme omega=k}{Let  $\varrho\in\M_1$. Uniformly for $x\geqslant 3$ and integer $h\geqslant 1$, we have
$$\EE_{x,\varrho}\big(\mu_2\chi_h\big)\ll 1.\eqdef{EM2chi}$$}
 \nid By \eqref{Pars}, we plainly have
$$\EE_{x,\varrho}\big(\mu_2\chi_h\big)\ll{1\over h!\log x}\int_0^1\bigg(\sum_{p\leqslant x}{\varrho(p)\{1+\cos(\vartheta\log p)\}\over p}\bigg)^h\d\vartheta.$$
In view of \eqref{condirho}, the contribution to the above integral of the range $0\leqslant \vartheta\leqslant 1/\log x$  is, for a suitable constant $C$ depending on $\varrho$, 
$$\ll {(2\log_2x+C)^h\over h!(\log x)^2}\ll1.$$ 
The complementary contribution is
$$\ll {1\over h!\log x}\int_{1/\log x}^1\big(\log_2x+\log 1/\vartheta+C\big)^h\d\vartheta={\e^C\over h!}\int_{\log_2x+C}^{2\log_2x+C}u^h\e^{-u}\d u\ll1.$$
\vskip-9mm\qed
\bigskip
\goodbreak
\par We are now in a position to complete the proof of \ref{thDelta2}.\par 
Let $\varepsilon>0$. For any $z\geqslant 1$, the contribution of $k\geqslant 4(1+\varepsilon)\log_2x$ to \eqref{decD2} is
$$\leqslant  \EE_{x,\varrho}\Big(z^{\omega(n)-4(1+\varepsilon)\log_2x}Ê\Delta^2\Big)\asymp{S_{2,\varrho_2}(x)\over  x(\log x)^{4(1+\varepsilon)\log z}}\ll (\log x)^{4z-3-4(1+\varepsilon)\log z}(\log_2x)^{ 4},
$$   with now $\varrho_2(n)=\varrho(n)z^{\omega(n)}$.
 Selecting $z=1+\varepsilon$, we get, for a suitable constant $c>0$,
$$\sum_{k\geqslant 4(1+\varepsilon)\log_2x} \EE_{x,\varrho}\big(\Delta^2\chi_k\big)\ll (\log x)^{1-c\varepsilon^2}(\log_2x)^{ 4}\ll\log x,$$
 provided $\varepsilon=\varepsilon_x:=\sqrt{5\log_3x/c\log_2x}$. 
\par
Similarly, bounding $\EE_{x,\varrho}(\Delta^2)$ with $z=1-\varepsilon_x$, we see that the contribution to \eqref{decD2}  of $k\leqslant 4(1-\varepsilon)\log_2x$ is $\ll\log x$.
\smallskip
Only $O\big(\sqrt{(\log_2x)\log_3x}\big)$ possible values of $k$ remain.
Each squarefree integer $n$ such that $\omega(n)=k$ may be represented as $n=ab$ with $\omega(a)=\lfloor k/2\rfloor$ and $\omega(b)=\plaf{ k/2}$.
However, since there are $\gg 2^{\omega(n)}/\sqrt{\log_2x}$ representations of $n$ as a product $ab$ of the above form, we infer that
$$\eqalign{\EE_{x,\varrho}(\Delta^2\chi_k)&\ll (\log x) \sqrt{\log_2x}\,\EE_{x,\varrho}\big(\mu_2\chi_{\fl{k/2}}\big)
\EE_{x,\varrho}\big(\mu_2\chi_{\plaf{ k/2}} \big)\ll (\log x) \sqrt{\log_2x},}$$
where the last estimate follows from \eqref{EM2chi}.
Summing over $k$, we get \eqref{S21x} as required. 
 
 \par   \smallskip\goodbreak   \bigskip  
\bigskip
 \centerline{\twelvebf References}
\bigskip
 {\eightpoint\leftskip9mm\rightskip5mm
\bibtem{BT12} R. de la Bretche \& G. Tenenbaum, Oscillations localisŽes sur les diviseurs, {\it J. London Math. Soc. (2) \bf 85 \rm (2012)}, 669--693.\par 
 \bibtem{BT23} R. de la Bretche \& G. Tenenbaum, Two upper bounds for the Erd\H os--Hooley Delta-function, {\it Sci. China Math. \bf66}, no. 12 (2023), 2683--2692.\par 
 \bibtem{BT25} R. de la Bretche \& G. Tenenbaum, Note on the mean-value of the Erd\H os-Hooley Delta-function, {\it Acta Arith.}, {\bf 219.4} (2025), 379--394.\par 
 \bibtem{BT26} R. de la Bretche \& G. Tenenbaum, Mean values of arithmetical functions and application to sum of powers, {\it Math. Proc. Cambridge Philos. Soc.}, {\bf 180}, no. 1 (2026), 1--13.\par
 \bibtem{Er73} P. Erd\H os, Problem 218, Solution by the proposer, {\it Can. Math. Bull. \bf 16} (1973),  621--622.\par 
 \bibtem{EN75} P. Erd\H os \& J.-L. Nicolas, RŽpartition des nombres superabondants, {\it Bull. Soc. math. France \bf 103} (1975), 65--90.\par 
 \bibtem{FGK23} K. Ford,  B. Green, D. Koukoulopoulos, Equal sums in random sets and the concentration of divisors, 
{\it Inventiones Math.} {\bf 232}, (2023),1027--1160. \par 
\bibtem{FKT23}  K. Ford, D. Koukoulopoulos, T. Tao, A lower bound on the mean value of the Erd\H os--Hooley Delta function, {\it Proc. London Math. Soc.} (3) {\bf 129, \rm no. 1 (2024), e12618.}
\par 
 \bibtem{HT86} R.R. Hall \& G. Tenenbaum, The average orders of Hooley's $\Delta_r$-functions, II,
 {\it Compos. Math. \bf 60} (1986), 163--186.\par 
 \bibtem{HT88} R.R. Hall \& G. Tenenbaum, {\it Divisors}, Cambridge tracts in
mathematics, no 90, Cambridge University Press (1988).\par 
\bibtem{Ho79} C. Hooley, On a new technique and its applications to the theory of numbers, {\it Proc.
London Math. Soc.} (3) {\bf 38} (1979), 115--151.
\par
  \bibtem{KT23} D. Koukoulopoulos \& T. Tao, An upper bound on the mean value of the Erd\H os--Hooley Delta function, {\it Proc. London Math. Soc. (3) \bf 127}, no.6 (2023),  1865--1885.
  \bibtem{NT98} M. Nair \& G. Tenenbaum, Short sums of certain arithmetic functions, {\it Acta Math. \bf180},
(1998), 119--144.\par
 \bibtem{No76} K. Norton, On the number of restricted prime factors of an integer, I, {\it
Illinois J. Math. \bf20}, no. 4 (1976),  681--705.\par 
\bibtem{Te85} G. Tenenbaum,
Sur la concentration moyenne des diviseurs, {\it Comment. Math. Helvetici} {\bf 60} (1985), 411-428. 
\bibtem{Te86} G. Tenenbaum,
Fonctions $\Delta$ de Hooley et applications, SŽminaire de thŽorie des nombres, Paris 1984-85, Prog. Math. 63 (1986), 225--239.
 \par
\bibtem{Te90} G. Tenenbaum, Sur une question d'Erd\H os et Schinzel, II, {\it Invent. Math.} {\bf 99}
(1990), 215--224.
\bibtem{Te15} G. Tenenbaum, {\it Introduction to analytic and probabilistic number theory}, 3rd ed., Graduate Studies in Mathematics 163, Amer. Math. Soc. (2015). See also {\it Introduction ˆ la thŽorie analytique et probabiliste des
nombres}, 5th ed., Dunod, Sciences Sup (Paris), 2022.\par
\par }
{\leftskip2mm\rightskip-2cm\sevenrm
\gutter=4cm \doublecolumns
 \obeylines \baselineskip=7pt
RŽgis de la Bretche
UniversitŽ Paris CitŽ, Sorbonne UniversitŽ, CNRS,
Institut Universitaire de France,
Institut de Math. de Jussieu-Paris Rive Gauche
 F-75013 Paris  
France
\smallskip
{\seventt regis.de-la-breteche@imj-prg.fr}
\goodbreak
 GŽrald Tenenbaum
Institut \'Elie Cartan
Universit\'e de Lorraine
 BP 70239\par
54506 Vand\oe uvre-ls-Nancy Cedex
 France
\smallskip
{\seventt gerald.tenenbaum@univ-lorraine.fr}
\singlecolumn
\par}
\end